\documentclass[12pt]{amsart} \usepackage{amscd} \usepackage{amsmath}

\newtheorem{theorem}{Theorem}
\newtheorem{lemma}[theorem]{Lemma}
\newtheorem{corollary}{Corollary}
\newtheorem{proposition}[theorem]{Proposition}

\theoremstyle{definition}

\theoremstyle{definition} \newtheorem{remark}{Remark}

\usepackage{amsfonts}
\newcommand{\field}[1]{\mathbb{#1}}          \newcommand{\Q}{\field{Q}}
\newcommand{\R}{\field{R}}                   \newcommand{\Z}{\field{Z}}
\newcommand{\C}{\field{C}}

 \renewcommand{\t}{\tau}

 \newcommand{\z}{\zeta}

\newcommand{\fg}{\mathfrak     g}     
     \newcommand{\fh}{\mathfrak    h}
\newcommand{\fu}{\mathfrak u} \newcommand {\ft}{\mathfrak t}

 \newcommand{\ra}{\rightarrow}

\begin{document}

\title[Generators  for  Arithmetic  Groups]{Generators for  Arithmetic
Groups}

\author{R. Sharma and T. N. Venkataramana}

\address{School   of  Mathematics,   Tata  Institute   of  Fundamental
Research, Homi Bhabha Road, Colaba, Mumbai 400005, India}

\email{ritumony@math.tifr.res.in, venky@math.tifr.res.in}

\subjclass{Primary 20F05, 11F06; Secondary 22E40}

\date{}

\begin{abstract}

We prove that  any non-cocompact irreducible lattice in  a higher rank
real semi-simple Lie  group contains a subgroup of  finite index which
is generated by {\bf three} elements.

\end{abstract}

\maketitle

\section{Introduction}

In  this paper  we study  the  question of  giving a  small number  of
generators for  an arithmetic  group.  Our main  theorem says  that if
$\Gamma $  is a  higher rank arithmetic  group, which  is non-uniform,
then $\Gamma  $ has a  finite index subgroup  which has at  most THREE
generators.   There are  reasons to  believe that  the bound  three is
sharp.\\

Our proof makes use of  the methods and results of \cite{T},\cite{R 4}
on certain unipotent generators for non-uniform arithmetic higher rank
groups  as also the  classification of  absolutely simple  groups over
number fields.  The  question of a small number  of generators is also
motivated by  the congruence subgroup  problem (abbreviated to  CSP in
the sequel). We prove the  following theorem, the main theorem of this
paper.\\

\begin{theorem}\label{mainth} Every higher rank {\bf non-uniform} 
arithmetic group  $\Gamma$ has a  subgroup $\Gamma '$ of  finite index
which is generated by at the most {\bf three} elements.
\end{theorem}

The proof exploits the existence  of certain unipotent elements in the
arithmetic group. The  higher rank assumption ensures that  if $U$ and
$U^-$ are opposing unipotent  radicals of maximal parabolic subgroups,
and  $M$  is  their  common  normaliser,  then  $M(\Z)$  will  have  a
``sufficiently  generic''  semi-simple  element.   There  are  generic
elements in $U^{\pm}(\Z)$ which  together with this generic element in
$M(\Z)$ will  be shown  to generate, {\bf  in general},  an arithmetic
group.   This is  already the  case  for the  group $SL(2,O_K)$  where
$K/\Q$ is a  non-CM extension of degree greater  than one (see section
2). \\

If $\gamma  $ is the above  ``generic'' element, and  $u^+\in U^+$ and
$u^-\in U^-$  are also  ``generic'', then let  $\Gamma $ be  the group
generated by the n-th powers $\gamma ^n$, $~(u^+)^n$ and $(u^-)^n$ for
some integer $n$.  Clearly, $\Gamma  $ is generated by three elements.
It is easy to show that  any arithmetic subgroup of $G(\Z)$ contains a
group  of the form  $\Gamma $  for some  integer $n$.   The genericity
assumption  will be shown  to imply  that for  {\bf most  groups} $G$,
$\Gamma $  intersects $U^+(\Z)$ and  $U^-(\Z)$ in subgroups  of finite
index.  Then a Theorem of Tits (\cite{T}) for Chevalley Groups and its
generalisation to  other groups of  $\Q$-rank $\geq 2$  by Raghunathan
\cite{R 4} (see also \cite{V}  for the case when $\Q$-rank ($G$)=$1$),
implies that $\Gamma $ is of finite index in $G(\Z)$.\\

The proof  that $\Gamma  $ intersects $U^{\pm}(\Z)$  in a  lattice for
most groups,  is reduced to  the existence of  a torus in  the Zariski
closure  of $M(\Z)$  (the  latter group  is  not equal  to $M$)  whose
eigen-spaces  (with  a given  eigen-value)  on  the  Lie algebra  $Lie
(U^{\pm})$ are  one dimensional.   The existence of  such a  torus for
groups of $\Q-rank\geq 2$ is proved by a case by case check, using the
Tits diagrams (classification) of  simple algebraic groups over number
fields.   It turns out  that in  the case  of exceptional  groups, the
existence  of such  a torus  is ensured  by the  results  of Langlands
\cite{L} and Shahidi \cite{Sh} who (in the course of their work on the
analytic continuation  of certain intertwining  operators) analyse the
action of the  Levi subgroup $L$ on the Lie  algebra $Lie(U^+)$ of the
unipotent radical. \\

However, this approach  fails in many groups of  $\Q$-rank one or two;
in these  cases, we  will have to  examine the individual  cases (i.e.
their  Tits   diagram),  to  produce  an  explicit   system  of  three
generators.  Thus,  a large part of  the proof (and a  sizable part of
the paper), involves, in low rank groups, a case by case consideration
of the Tits  diagrams. In many of these cases,  the explicit system of
generators is  quite different from  the general case (see  sections 4
and 5). \\

We  end this  introduction with  some notation.   Given  a $\Q$-simple
semi-simple  algebraic group,  there  is an  absolutely almost  simple
algebraic  group ${\mathcal  G}$ over  a  number field  $K$ such  that
$G=R_{K/\Q}({\mathcal G})$ where $R_{K/\Q}$ is the Weil restriction of
scalars.  Moreover, $\Q-rank (G)=K-rank ({\mathcal G})$ and $G(\Z)$ is
commensurate to  $G(O_K)$ where $O_K$ is  the ring of  integers in the
number field. For  these reasons, we replace henceforth  the group $G$
over $\Q$  with an  absolutely simple group  (still denoted $G$  by an
abuse of notation), defined over a number field $K$.\\

Given a group  $G$, and element $g,h\in G$ and  a subset $S\subset G$,
denote by  $^g(h)$ the  conjugate $ghg^{-1}$, and  $^g(S)$ the  set of
elements $ghg^{-1}$ with $h\in S$. \\

If $\Gamma _0$  is a group, $\Gamma , \Delta$  are subgroups, one says
that $\Gamma  $ {\bf virtually  contains} $\Delta$ and  writes $\Gamma
\geq \Delta$ if the intersection $\Gamma \cap \Delta$ has finite index
in $\Delta$. One says that  $\Gamma $ is {\bf commensurate} to $\Delta
$ and  writes $\Gamma \simeq  \Delta$ if $\Gamma $  virtually contains
$\Delta$ (.i.e.   $\Gamma \geq \Delta$) and vice  versa (i.e.  $\Delta
\geq \Gamma $). \\

\section{Preliminary results on Rank one groups}

\subsection{The Group SL(2)}
In  this  subsection  we  prove  Theorem  \ref{mainth}  for  the  case
$G=EL(2)$  over a  number field  $E$.  The  assumption of  higher rank
translates  into the  condition that  $E$ has  infinitely  many units.
That is, $E$  is neither $\Q$ nor an  imaginary quadratic extension of
$\Q$.  It turns out that if $E$ is not a CM field, that is, $E$ is not
a  totally imaginary  quadratic  extension of  a  totally real  number
field, then, the  proof is easier.  We will  therefore prove this part
of the theorem first.

\begin{proposition}\label{noncm} Let $K$ be a number field, which is 
not $\Q$ and  which is not a CM  field. Let $G=R_{K/\Q}(EL(2))$. Then,
any  arithmetic subgroup  of $G(\Q)$  has a  subgroup of  finite index
which has three generators.
\end{proposition}

Before we begin the proof of  Proposition \ref{noncm} , we prove a few
Lemmata.  We will first assume that  $E$ is a non-CM number field with
infinitely many  units. Let $O_E$ denote  the ring of  integers in the
number field and  $O_E^*$ denote the multiplicative group  of units in
the ring $O_E$.

\begin{lemma}\label{finite} Let $\Delta $ be a subgroup of finite
index in $O_E^*$ and $F$  the number field generated by $\Delta$. Then
$F=E$.
\end{lemma}
\begin{proof} If $r_1(K)$ and $r_2(K)$ are the number of inequivalent 
real  and complex  embeddings  of a  number  field $K$,  then, by  the
Dirichlet Unit Theorem, the  rank of $O_E^*$ is $r_1(E)+r_2(E)-1$. Let
$d$ be the degree of $E$ over $F$. \\

Let $A$ be the set of real places of $F$. To each $a\in A$, let $x(a)$
be the  number of real  places of $E$  lying above $a$ and  $y(a)$ the
number of non-conjugate  complex places of $E$ lying  above $a$. Then,
for  each $a\in  A$ we  have $x(a)+2y(a)=d$,  the degree  of  $E$ over
$F$. Clearly, $x(a)+y(a)\geq 1$ for each $a$. \\

Let $B$ be the number of non-conjugate complex places of $F$. Then all
places of  $E$ lying above a  place $b\in B$ are  imaginary.  If their
number is $y(b)$ , then we have $y(b)=d$ for each $b$. \\

The rank of the group of units  $O_F^*$ of the number field $F$ is, by
the Drichlet Unit Theorem, $Card (A)+Card (B)-1$. That of $O_E^*$ is
\[-1+ \sum _{a\in A} (x(a)+y(a))+ \sum_{b\in B} y(b) \]
By assumption, $O_F^*$  and $O_E^*$ have the same  rank, since $O_F^*$
contains $\Delta$, a subgroup of finite index in $O_E^*$. We thus have
the equation
\[\label{card}(\ref{card})
~~~~Card(A)+Card(B)=\sum  _{a\in  A}(x(a)+y(a))+\sum_{b\in  B}y(b)  \]
Since $x(a)+y(a)\geq  1$ and $y(b)=d\geq 1$,  equation \ref{card} show
that if $B$ is non-empty, then $d=1$ and $E=F$. \\

If $B$ is empty, then $F$ has no complex places, and so $F$ is totally
real.   Moreover, since $x(a)+y(a)\geq  1$, equation  \ref{card} shows
that  for each  $a\in  A$, $x(a)+y(a)=1$.   Thus,  either $x(a)=0$  or
$y(a)=0$.  If, for some $a$, $y(a)=0$ then the equation $d=x(a)+2y(a)$
shows that $d=1$ and $E=F$. \\

The only possibility left is that $x(a)=0$ and $y(a)=1$ for each $a\in
A$,  and  $F$  is  totally  real.   Therefore,  for  each  archimedean
(necessarily real) place $a$ of $F$, we have $y(a)=1$ and $d=2y(a)=2$,
that is there  is only one place  of $E$ lying above the  place $a$ of
$F$ and  is a  complex place, whence  $E/F$ is a  quadratic extension,
which is  totally imaginary.   Hence $E$ is  a CM-extension,  which is
ruled out by assumption.
\end{proof}

The field  extension $E$ over $\Q$  has only finitely  many proper sub
fields  $E_1, E_2,\cdots,  E_m$  (this follows  trivially from  Galois
Theory, for example).\\

\begin{lemma} \label{exist} Suppose that $E$ is a number field which 
is not  a CM field.  Then  There exists an element  $\theta \in O_E^*$
such that for  any integer $r\geq 1$, the sub  ring $\Z[\theta ^r]$ of
$O_E$ generated  by $\theta ^r$ is  a subgroup of finite  index in the
additive group  $O_E$. In particular, $\Z[\theta  ^r]\supset NO_E$ for
some integer  $N$. Consequently, there  exists an element  $\theta \in
O_E^*$ which does not lie  in any of the subfields $E_1,\cdots,E_m$ as
above, and for every such $\theta $, the sub-ring $\Z[\theta ^r]$ is a
subgroup of of finite index in $O_E$.
\end{lemma}
\begin{proof} By Lemma \ref{finite} the intersection 
$\Delta _i=O_E^*\cap E_i$ is of  {\bf infinite index} in $O_E^*$.  Let
us now write the abelian  group $O_E^*$ additively.  Then, we have the
$\Q$  -subspaces  $W_i=\Q\otimes  \Delta   _i$  of  the  vector  space
$W=\Q\otimes  O_E^*$ (the latter  of dimension  $r_1(E)+r_2(E)-1$ over
$\Q$).  Since $W_i$ are proper subspaces of $W$, it follows that there
exists an element  of $W$ (hence of the  subgroup $O_E^*$) no rational
multiple  of which  lies  in  $W_i$ for  any  $i$.  Interpreting  this
statement  multiplicatively,  there exists  an  element  $\theta $  of
$O_E^*$  such that  no integral  power of  $\theta $  lies in  the sub
fields $E_i$  for any $i$.   Consequently, for any integer  $r\neq 0$,
the subfield  $\Q[\theta ^{r\Z}]$ is  all of $E$.  In  particular, the
sub ring $\Z[\theta  ^r]$ generated by $\theta ^r$  is of finite index
in the ring $O_E$.
\end{proof}

We  now  begin the  proof  of  Proposition  \ref{noncm}. Consider  the
matrices  $u_+=\begin{pmatrix}   1  &   1  \\  0   &  1\end{pmatrix}$,
$u_-=\begin{pmatrix} 1 &  0 \\ 1 & 1\end{pmatrix}$ and  by an abuse of
and by an  abuse of notation, $\theta =u_+=\begin{pmatrix}  \theta & 0
\\ 0 &  \theta ^{-1}\end{pmatrix}$, where $\theta \in  O_K^*$ is as in
Lemma \ref{exist}.  Then,  the group $\Gamma =<u_+^r,u_-^r,\theta ^r>$
generated  by  $u_{\pm}^r$  and  $\theta ^r$  contains,  for  integers
$m_1,m_2,\cdots m_l$, and $n_1,n_2,\cdots n_l$, the element
\[^{\theta ^{m_1r}}(u_+^{rn_1}) ^{\theta ^{m_2r}}(u_+^{rn_2})
\cdots ^{\theta  ^{m_lr}}(u_+^{rn_l}).  \] This element  is simply the
matrix   $\begin{pmatrix}  1   &  r\sum   n_i\theta  ^{2m_ir}\\   0  &
1\end{pmatrix}$.   Picking  suitable $m_i,  n_i$,  we  get from  Lemma
\ref{exist}, an  integer $N$ such that  $\begin{pmatrix} 1 & x  \\ 0 &
1\end{pmatrix}\in \Gamma $ for  all $x\in NO_K$.  Similarly, all lower
triangular  matrices  of  the form  $\begin{pmatrix}  1  &  0 \\  x  &
1\end{pmatrix}$  are in $\Gamma  $ for  all $x\in  NO_K$. But  the two
subgroups $U^+(NO_K)=\begin{pmatrix}  1 & NO_K \\  0 & 1\end{pmatrix}$
and $U^-(NO_K)=\begin{pmatrix}  1 & 0 \\  NO_K & 1\end{pmatrix}\subset
\Gamma  $ generate  a subgroup  of finite  index in  $SL_2(O_K)$  ( by
\cite{Va}). Hence $\Gamma  $ is of finite index  in $SL_2(O_K)$. It is
clear  that any  subgroup of  finite index  in $SL_2(O_K)$  contains a
three  generated group  $\Gamma  =<u_+^r,\theta ^r,  u_-^r>$ for  some
$r$. This completes the proof of Proposition \ref{noncm}.

\subsection{The CM case.} Suppose that $F$ is a totally real number
field  of  degree  $k\geq 2$  and  suppose  that  $E/F$ is  a  totally
imaginary quadratic extension of  $F$. There exists an element $\alpha
\in E$ such that $\alpha ^2=-\beta  \in F$ where $\beta $ is a totally
positive element  of $F$  (that is,  $\beta $ is  positive in  all the
archimedean  (hence real)  embeddings of  $F$).  Let  $\theta $  be an
element of infinite order in  $O_F^*$ as in Lemma \ref{exist}, so that
for  any integer  $r$,  the sub-ring  $\Z[\theta  ^r]$ of  $O_F$ is  a
subgroup of finite  index in $O_F$ (in Lemma  \ref{exist}, replace $E$
by the totally real field $F$). We have thus the following analogue of
Lemma \ref{exist}, in the CM case.
\begin{lemma} \label{exist'} 
Suppose that  $E$ is a CM  field and is a  totally imaginary quadratic
extension of a totally real  number field $F$. There exists an element
$\theta  \in O_E^*$  such that  for any  integer $r\neq  0$,  the ring
$\Z[\theta ^r]$ generated by $\theta ^r$ is a subgroup of finite index
in $O_F$.
\end{lemma}

\begin{proof} The group of units of $E$ contains the group of units of
$F$  as a  subgroup  of finite  index.   Therefore, we  may apply  the
previous  lemma (Lemma  \ref{exist}), with  $E$ replaced  by  $F$ (the
latter is not a CM field).
\end{proof} 

Consider the elements $h=h(\theta)=
\begin{pmatrix}\theta & 0 \\ 0 & \theta ^{-1}\end{pmatrix}$, $u_+=
\begin{pmatrix}1 & 1 \\ 0 & 1\end{pmatrix} $, and $u_-=
\begin{pmatrix}1 & 0 \\ \alpha & 1\end{pmatrix}$ of $SL(2,O_E)$. 
Given an arithmetic subgroup  $\Gamma _0$ of $SL(2,O_E)$, there exists
an  integer $r$  such that  the  group $\Gamma  =<h^r, u_+^r,  u_-^r>$
generated  by the  $r$-th  powers  $h^r, u_+^r$  and  $u_-^r$ lies  in
$\Gamma _0$.

\begin{proposition}\label{CM}For every integer $r$, the group 
$\Gamma $ in the foregoing paragraph is arithmetic (i.e.  is of finite
index in  $SL(2,O_E)$).  In  particular, every arithmetic  subgroup of
$SL(2,O_E)$ is virtually $3$-generated.
\end{proposition}

\begin{proof} Write the Bruhat decomposition for the element 
\[u_-^r=\begin{pmatrix}1 & 0\\ r\alpha & 1\end{pmatrix}=
\begin{pmatrix}1 & \frac{1}{r\alpha} \\ 0 & 1\end{pmatrix}
\begin{pmatrix}-\frac{1}{r\alpha} & 0 \\ 0 & -r\alpha\end{pmatrix}
\begin{pmatrix}0 & 1 \\ -1 & 0\end{pmatrix}
\begin{pmatrix}1 & \frac{1}{r\alpha} \\ 0 & 1\end{pmatrix}\]

By  the  choice of  the  element $\theta  $,  the  group generated  by
$h(\theta )^r$ and $u_+^r$ contains the subgroup $U^+(NO_F)=\{u=
\begin{pmatrix}1 & Nb \\ 0 & 1\end{pmatrix}: b\in O_F\}$.  Clearly,
$\Gamma  \supset  U^+(O_F)$.   Define  $U^-(NO_F)$  similarly.   Since
$\alpha ^2$ lies in the  smaller field $O_F$, a computation shows that
the    conjugate   $^{u_-^r}(U^+(NO_F))$    contains    the   subgroup
$^{v_+}(U^-(N'O_F))$ for some integer $N'$, where $v_+$ is the element
$\begin{pmatrix}1  & \frac{1}{r\alpha} \\  0 &  1\end{pmatrix}$.  Thus
the group $^{v_+}(U^-(N'O_F))\subset \Gamma  $.  Since the group $U^+$
is  commutative, we have  $^{v_+}(U^+(N'O_F))=U^+(N'O_F)\subset \Gamma
$.  By  a Theorem  of Vaserstein (\cite{Va}),  the group  generated by
$U^+(NO_F)$  and  $U^-(N'O_F)$  is  a  subgroup  of  finite  index  in
$SL(2,O_F)$. In particular, it  contains some power $h^M=h(\theta )^M$
of $h$.\\

Since a power of  $h$ already lies in $\Gamma $, we  see that for some
integer $r'$,  the commutator  $u_1$ of $h^{r'}$  and $^{v_+}(h^{r'})$
lies  in  $\Gamma  $.   This  commutator is  nothing  but  the  matrix
$\begin{pmatrix}1  &  (\theta   ^{2r'}-1)(\frac{1}{r\alpha})  \\  0  &
1\end{pmatrix}$.   Now,  $\frac  {1}{\alpha  }=\frac{-\alpha}{\beta}$,
with  $\beta \in F$.   Therefore, by  Lemma \ref{exist},  the subgroup
generated by  $\theta ^{r'}$  and $u_1$ contains,  for some  $M'$, the
subgroup  $U^+(M'O_F\alpha)$  consisting   of  elements  of  the  form
$\begin{pmatrix}1  &   xM'\alpha  \\  0   &  1\end{pmatrix}$.   Hence,
$U^+(M'O_F\alpha  )\subset \Gamma  $.   We have  already  seen in  the
beginning of the last  paragraph that $U^+(MO_F)\subset \Gamma $. Now,
up  to a  subgroup of  finite index,  $O_E$ is  the sum  of  $O_F$ and
$O_F\alpha$  since  $E/F$  is   a  quadratic  extension  generated  by
$\alpha$. This shows  (after changing $M'$ if necessary  by a suitable
multiple), that  $U^+(M'O_E)\subset \Gamma $.  The  conjugate of $U^+$
by the {\it lower triangular}  matrix $\begin{pmatrix}1 & 0\\ \alpha &
1\end{pmatrix}$  is  a  unipotent  group  $V$ opposed  to  $U^+$.   By
Vaserstein's  Theorem in  \cite{Va} (replacing  $U^-$ by  our opposite
unipotent  group $V$)  we see  that $\Gamma  $ contains  an arithmetic
subgroup of $SL(2,O_E)$ and is hence itself arithmetic.
\end{proof}

For handling  some $\Q$-rank one groups,  we will need  a more general
version  of  Proposition \ref{CM}.   Let  $x\in  E\setminus  F$ be  an
integral  element divisible  by $N!$  (the  product of  the first  $N$
integers) for a large rational integer $N$.  Denote by $U^-(xO_F)$ the
set of matrices of the form $\begin{pmatrix} 1 & 0\\ xa & 0
\end{pmatrix}$ with $a\in O_F$.  Denote by $U^+(rO_F)$ the group of
matrices of the form $\begin{pmatrix} 1 & ra\\ 0 & 1
\end{pmatrix}$ with $a\in O_F$. 

\begin{proposition}\label{CM'}
The group generated by $U^+(rO_F)$  and $U^-(xO_F)$ is of finite index
in $SL(2,O_E)$.
\end{proposition} 

\begin{proof} Denote by $\Gamma $ the group in the proposition. 
We first find an element  $\begin{pmatrix}a & b\\ c & d \end{pmatrix}$
in $\Gamma  $ such that $ac\neq 0$  and $c$ lies in  the smaller field
$F$. To  do this,  we use  the existence of  infinitely many  units in
$F$. Write $x ^2=tx -n$ with $t(=tr_{E/F}(x))$ and $n(=N_{E/F}(x))$ in
$F$. Assume  that $t\neq  0$ since $t=0$  has already been  covered in
Proposition  \ref{CM}.  Given  a unit  $\theta\in O_E^*$  consider the
product element $g\in SL(2,E)$ given by $g=
\begin{pmatrix} 1 & 0 \\ -x\theta ^{-1} & 1\end{pmatrix} 
\begin{pmatrix} 1 & \frac{\theta -1}{t}\\ 0 & 1 \end{pmatrix} 
\begin{pmatrix} 1 & 0 \\ x & 1\end{pmatrix}$. \\

Now,  normal subgroups  of  higher rank  arithmetic  groups are  again
arithmetic    (\cite{M},    \cite{R    1},   \cite{R    2}).     Since
$T(O_F)(=T(O_E))$ normalises  the groups $U^+(rO_F)$  and $U^-(xO_F)$,
it follows that to prove the  arithmeticity of $\Gamma $, it is enough
to prove  the arithmeticity  of the group  generated by $\Gamma  $ and
$T(O_F)$.  We may thus assume that $\Gamma $ contains $T(O_E)=T(O_F)$.
Here the equality is up to subgroups of finite index. \\

If $\theta $ is a unit  such that $\theta \equiv 1$ (mod $tO_F$), then
from the definition of $g$ and $\Gamma $ it is clear that $g\in \Gamma
$.   Write  $g=  \begin{pmatrix}  a  & b\\  c  &  d\end{pmatrix}$.   A
computation  shows that $a=1+\frac{\theta  -1}{t}x$, $c=\frac{1-\theta
^{-1}}{t}n$.   Since $x$  and $1$  are linearly  independent  over $F$
($x\notin F$), it follows that $a\neq 0$ and in fact that $a\notin F$.
The expression for $c$ shows that $c\neq 0$.\\

The  Bruhat   decomposition  for  $g=\begin{pmatrix}  a  &   b\\  c  &
d\end{pmatrix}$ is given by
\[g=\begin{pmatrix} a & b\\ c & d\end{pmatrix}=
\begin{pmatrix} 1 & ac^{-1}\\ 0 & 1\end{pmatrix} 
\begin{pmatrix} c^{-1} & 0\\ 0 & c\end{pmatrix}
\begin{pmatrix} 0 & -1\\ 1 & 0\end{pmatrix}
\begin{pmatrix} 1 & dc^{-1}\\ 0 & 1\end{pmatrix}. \] 

Thus, $\Gamma \supset {}^{\big (\begin{smallmatrix}a & b\\ c & d
\end{smallmatrix}\big )}
\begin{pmatrix} 1 & rO_F\\ 0 & 1\end{pmatrix}= 
^{\big  (\begin{smallmatrix}1 &  ac^{-1}\\ 0  & 1\end{smallmatrix}\big
)}\begin{pmatrix} 1 & 0 \\c^2rO_F & 1\end{pmatrix}$.  Moreover $\Gamma
\supset
\begin{pmatrix} 1 & rO_F\\ 0 & 1\end{pmatrix}=
^{\big (\begin{smallmatrix}1  & ac^{-1}\\ 0  & 1 \end{smallmatrix}\big
)}
\begin{pmatrix} 1 & rO_F\\ 0 & 1\end{pmatrix}$. The group $\Delta $ 
generated  by $\begin{pmatrix}  1  & rO_F\\  0  & 1\end{pmatrix}$  and
$\begin{pmatrix} 1 & 0\\  c^2rO_F & 1\end{pmatrix}$ contains, for some
integer $r'$, $U^+(r'O_F)$ and $U^-(r'O_F)$ (since $c$ and hence $c^2$
lie  in the  field  $F$) Hence,  by  Vaserstein's Theorem  (\cite{Va})
$\Delta $ is an arithmetic subgroup of $SL(2,O_F)$. \\

In   particular,    $\Gamma   $   contains    the   subgroup   $^{\big
(\begin{smallmatrix}  1  &  ac^{-1}\\  0 &  1\end{smallmatrix}\big  )}
(\theta  ^{r''\Z})$ for  some integer  $r''$.  By  enlarging  $r''$ if
necessary, assume that $\theta  ^{r''\Z}\subset \Gamma $. Thus $\Gamma
$ contains the commutator group
\[ [\begin{pmatrix} 1 & ac^{-1}\\ 0 & 1\end{pmatrix}, \theta
^{r''\Z}]= \begin{pmatrix} 1 & ac^{-1}(\sum \Z(\theta ^{r''k}-1)) \\ 0
& 1\end{pmatrix}\]  where the  sum is over  all integers $k$.   By the
properties of  the element $\theta$, the  sum is a  subgroup of finite
index in  the ring $O_F$,  whence, we get  an integer $r_0$  such that
$\Gamma   \supset   \begin{pmatrix}  1   &   ac^{-1}r_0   O_F\\  0   &
1\end{pmatrix}$. Since $\Gamma $ already contains $\begin{pmatrix} 1 &
rO_F\\  0  & 1\end{pmatrix}$,  $a$  does not  lie  in  $F$, and  $O_E$
contains $O_F\oplus  r_0ac^{-1}O_F$ for  a suitable $r_0$,  it follows
that for a  suitable integer $r_1$, the subgroup  $\begin{pmatrix} 1 &
r_1O_E\\  0 &  1\end{pmatrix}$ lies  in $\Gamma  $. Now  $\Gamma  $ is
obviously  Zariski dense  in $SL(2,O_E)$;  moreover it  intersects the
unipotent radical  $U^+$ is an  arithmetic group. Hence  it intersects
some opposite unipotent  radical also in an arithmetic  group; but two
such opposing unipotent arithmetic groups generate an arithmetic group
(\cite {Va}). Therefore, $\Gamma $ is arithmetic.
\end{proof}

\subsection{the group SU(2,1)}
In this section, we prove results on the group $SU(2,1)$ (with respect
to a  quadratic extension $L/K$ of  a number field  $K$), analogous to
those in the  section on $SL_2$. These will be needed  in the proof of
Theorem 1, in those cases where a suitable $SU(2,1)$ embeds in $G$. \\

Suppose   that   $E/\Q$  is   a   {\bf   real}  quadratic   extension,
$E=\Q(\sqrt{z})$ with  $z >0$.  Denote  by $x\mapsto \overline  x$ the
action of the non-trivial element  of the Galois group of $E/\Q$.  Let
$h=\begin{pmatrix}0 & 0 & 1\\0 & 1  & 0\\ 1 & 0 & 0\end{pmatrix}$.  We
will view $h$ as a form in three variables on $E^3$ which is hermitian
with respect to this non-trivial Galois automorphism.  Set
\[G=SU(h)=SU(2,1)=\{g\in SL_3(E):\overline {^t g}hg=h\}.\] 
Then $G$ is an algebraic group over $\Q$.\\

Define the groups
\[U^+=\{\begin{pmatrix}1 & z & -\frac{z\overline z}{2}\\ 0 & 1 &
-\overline z \\ 0 & 0 & 1 \end{pmatrix}
\begin{pmatrix}1 & 0 & w\\ 
0  &  1  &  0  \\  0  & 0  &  1  \end{pmatrix}:  w+\overline  w=0\},\]
$U^-=^t(U^+)$, the subgroup of $SU(2,1)$ which is an opposite of $U^+$
consisting of matrices which are  transposes of those in $U^+$ and let
$T$ be  the diagonals  in $SU(2,1)$. Then,  up to subgroups  of finite
index, we  have $T(\Z)=\{\begin{pmatrix}\theta  & 0 &  0\\ 0  & \theta
^{-2} & 0  \\ 0 & 0 & \theta \end{pmatrix}:  \theta \in O_E^*\}$. Note
that for a  unit $\theta \in O_E^*$, we  have $\theta \overline \theta
=\pm 1$. \\

Suppose that $F/\Q$ is imaginary quadratic, $t\in O_F\setminus \Z$ and
{\bf define} the group $U^+(t\Z)$ as the one generated by the matrices
$\begin{pmatrix}1 & 0 & tx\sqrt{z} \\ 0 & 1 & 0 \\ 0 & 0 & 1
\end{pmatrix}$, and $\begin{pmatrix}1 & tux & -\frac{t^2 x^2
u\overline u}{2}\\ 0 & 1 & -t\overline {u} x\\ 0 & 0 & 1
\end{pmatrix}$ with $x\in \Z$ and $u\in O_E$. Denote by $U_{2\alpha}$
the root group corresponding to the root $2\alpha$, where $\alpha $ is
the  simple  root  for  ${\bf  G}_m (\subset  T)$  occurring  in  $Lie
U^+$. Here,  the inclusion of ${\bf G}_m$  in $T$ is given  by the map
$x\mapsto  \begin{pmatrix}  x &  0  &  1\\ 0  &  1  & 0  \\  0  & 0  &
x^{-1}\end{pmatrix}$. Note  that the commutator  $[U^+(t\Z), U^+(\Z)]$
is  $U_{2\alpha }(t\Z)\subset  U^+(t\Z)$.  Hence  $U^+(\Z)$ normalises
$U^+(t\Z)$.   Note  moreover  that  $U^+(t\Z)$ contains  the  subgroup
$U_{2\alpha}(t^2  \Z+t\Z))$;  now,  the  elements $t$  and  $t^2$  are
linearly independent  over $\Q$, hence $t\Z+t^2\Z$  contains $r\Z$ for
some integer $r>0$. \\
 
\begin{proposition}If $\Gamma \subset G(O_F)$ is such that for some
$r\geq  1$,  the group  $\Gamma  $  contains  the group  generated  by
$U^+(rt\Z)$  and $U^-(r\Z)$,  then $\Gamma  $  is of  finite index  in
$G(O_F)$.
\end{proposition}

\begin{proof} By the last remark in the paragraph preceding the 
proposition, there exists an integer,  we denote it again by $r$, such
that   the   group   $\Gamma   $   contains   $U_{2\alpha}(r\Z)$   and
$U^-(r\Z)$. Thus, by \cite{V}  (note that $\R-rank (G)=2$ since $E/\Q$
is real quadratic and $G(\R)=SL(3,\R)$), $\Gamma $ contains a subgroup
of $SU(2,1)(\Z)$  of finite index.  Therefore, $\Gamma  $ contains the
group generated  by $U^+(rt\Z)$ and  $U^+(r\Z)$ for some  integer $r$.
Clearly  the group  generated contains  $U^+(r'O_F)$ for  some integer
$r'$  (since  $F/\Q$  is  quadratic  and  $t$  and  $1$  are  linearly
independent over $\Q$). Therefore,  by \cite{V} again, we get: $\Gamma
$ is an arithmetic subgroup of $G(O_F)$.
\end{proof}

We now prove a slightly stronger version of the foregoing proposition.
\begin{proposition} \label{SU(2,1)} 
Suppose  that   $E$  and  $F$   are  as  before,   $E=\Q\sqrt{z}$  and
$F=\Q(\sqrt{t})$. Let  $\Gamma \subset G(O_F)$  be such that  for some
integer   $r$,  $\Gamma   $   contains  the   groups  $U^-(r\Z)$   and
$U_{2\alpha}(rt\Z)$. Then, $\Gamma $ is of finite index in $G(O_F)$.
\end{proposition}

\begin{proof} Consider the map $f:SL(2)\ra SU(2,1)$ given by 
  $\begin{pmatrix} a & b\\ c & d\end{pmatrix}\mapsto \begin{pmatrix} a
&   0   &  b\sqrt{z}\\0   &   1   &   0\\\frac{c}{\sqrt{z}}  &   0   &
d\end{pmatrix}$. The  map $f$  is defined over  $\Q$, takes  the upper
triangular matrices with 1s on the diagonal to the group $U_{2\alpha}$
and takes the Weyl group element  $w$ into the $3\times 3$ matrix $w'$
which has  non-zero entries on the anti-diagonal  and zeros elsewhere.
Under conjugation action by the element $f(h)$ with $h=\begin{pmatrix}
a &  0\\0 & a^{-1}\end{pmatrix}$,  the group $U^+(r\Z)$ is  taken into
$U^+(ra\Z)$. Under  conjugation by $w'$,  $U^{-}$ is taken  into $U^+$
and vice versa. \\

Write  the Bruhat decomposition  of $u_+=\begin{pmatrix}  1 &  rt\\0 &
1\end{pmatrix}$, with  respect to the lower triangular  group.  We get
$u_+=v_1^{-}h_1wu_1^{-}$.   Here,  $h_1=\begin{pmatrix}-rt  &  0\\0  &
-\frac{1}{rt}\end{pmatrix}$. If $r$ is  suitably large, then $\Gamma $
contains $u_+$ by assumption.  To prove arithmeticity, we may assume (
see  the   proof  of  Proposition  \ref{CM'})   that  $\Gamma  \supset
T(r\Z)$. Then,
\[\Gamma \supset ^{u_+}(U^-(r\Z))\supset
^{v_1^-}(U_{\alpha}(rt\Z))\supset    ^{v_1^-}(U_{2\alpha}(r^2t^2\Z)).\]
The  last  inclusion follows  by  taking  commutators  of elements  of
$U_{\alpha}(rt\Z)$ where, $U_{\alpha}(rt\Z)$ is the group generated by
the  elements $\begin{pmatrix}  1 &  rt &  -\frac{r^2t^2}{2}\\0 &  1 &
-rt\\ 0 &  0 & 1 \end{pmatrix}$. Note that  $t^2\in \Q$ by assumption.
Hence $\Gamma  \supset ^{v_1^-}(U_{2\alpha }(r'\Z))$  for some integer
$r'$. \\

Define  $U_{-\alpha}(rt\Z)$  similarly  to  the above  (e.g.   as  the
transpose of $U_{\alpha}$). Since $v_1^-$ centralises all of $U^-$, we
obtain         $\Gamma         \supset         U_{-\alpha}(r\Z)\supset
^{v_1^-}(U_{-\alpha}(r\Z))$. \\
 
The conclusions  of the  last two paragraphs  and \cite{V}  shows that
there exists  a subgroup  $\Delta $ of  finite index  in $SU(2,1)(\Z)$
such that  $\Gamma \supset ^{v_1^-}(\Delta)$. In  particular, for some
integer $r'$,  the group $^{v_1^-}(\Gamma  )$ contains {\it  both} the
groups $U_{\alpha}(r't\Z)$  and $U_{\alpha }(r\Z)$.   Consequently, it
contains  $U^+(r'O_F)$, a  subgroup of  finite index  in  the integral
points of  the unipotent  radical of a  minimal parabolic  subgroup of
$SU(2,1)$  over  $F$ (note  that  up  to  subgroups of  finite  index,
$t\Z+\Z=O_F$). Note also that  the real rank of $SU(2,1)(F\otimes \R)$
is at least two.  Now, $\Gamma $ is clearly Zariski dense in $SU(2,1)$
regarded  as a  group over  $F$. Therefore,  by \cite{V},  $\Gamma$ is
arithmetic. \end{proof}

\subsection{Criteria for Groups of Rank One over Number Fields}
Suppose that  $K$ is a number  field. Let $G$ be  an absolutely almost
simple algebraic group with  $K$-rank ($G$)$\geq 1$. Let $S\simeq {\bf
G}_m$ be  a maximal $K$-split torus  in $G$, $P$  a parabolic subgroup
containing $S$, and $U^+$ the  unipotent radical of $P$. Let $M\subset
P$  be the  centraliser of  $S$ in  $G$.  Let  $M_0$ be  the connected
component  of identity  of the  Zariski  closure of  $M(O_K)$ in  $M$.
Write  $\fg$ for  the Lie  algebra  of $G$.   We have  the root  space
decomposition  $\fg= \fg  _{\pm \alpha}\oplus  \fg _0$,  where $\oplus
_{\alpha >0}\fg  _{\alpha}$ is a decomposition  of $ Lie  U^+$ for the
adjoint action of $S$. Denote by  $log :U^+\ra \fu$ the log mapping on
the unipotent group $U^+$. It  is an isomorphism of $K$-varieties (not
of  groups in  general). Define  similarly $U^-$  to be  the unipotent
$K$-group  group  with  Lie  algebra  $\fu  ^-=\oplus  {\alpha  >0}\fg
_{-\alpha}$. This is the ``opposite'' unipotent group. There exists an
element $w\in N(S)/Z(S)$  in the Weyl group of  $G$ ($N(S)$ and $Z(S)$
being  the  normaliser and  the  centraliser  of  $S$ in  $G$),  which
conjugates $U^+$ into $U^-$.  Further, the map $(u,m,v)\mapsto umwv=g$
maps $U^+\times M\times U^+$ isomorphically onto a Zariski open subset
of $G$. \\

The  following technical proposition  will be  used repeatedly  in the
sequel.
\begin{proposition} \label{technical} Suppose that $K$ is any number
field. Let $G$  be of $K$-rank $\geq 1$.   Let $\Gamma \subset G(O_K)$
be  Zariski dense, and  assume that  $\R$-rank ($G_\infty$)  $\geq 2$.
Suppose that there exists an element $m_0\in M(O_K)$ of infinite order
such that 1)  all its eigenvalues are of infinite  order in its action
on $Lie U^+$, 2) if $g=umwv\in \Gamma $ , then there exists an integer
$r\neq  0$ such  that  $^u(m_0^r)\in  \Gamma $.   Then,  $\Gamma $  is
arithmetic.
\end{proposition} 

\begin{proof} Let $V$ be the Zariski closure of the 
intersection of $U^+$  with $\Gamma$.  View $U^+$ as  a $\Q$-group, be
restriction of  scalars.  By  assumption, for a  Zariski dense  set of
elements $u\in  U^+(\Q)$, there exists  an integer $r=r(u)$  such that
the commutator  $[m_o^r,u]$ lies in  $\Gamma $.  If $\frak  v$ denotes
the $\Q$-Lie algebra of $V$, then, this means that, $\frak v$ contains
vectors of the  form $(Ad(m_0^r)-1)(log u)$ with $log  u$ spanning the
$\Q$-vector space  $\fu$. By  fixing finitely many  $u'$ which  give a
basis of $\fu$ (as a $\Q$-vector  space), we can find a common integer
$r$ such that $(Ad (m_0^r)-1)log u\in  \frak v$, for all $u$; in other
words, $(Ad(m_0^r)-1)(\fu)\subset  \frak v$.  The  assumption on $m_0$
now  implies that  $\frak  v=\fu$.  Hence  $V=U^+$,  which means  that
$\Gamma  \cap U^+\subset  U^+(O_K)$  is Zariski  dense  in $U^+$.   By
\cite{R 5}, Theorem  (2.1), it follows that $\Gamma  \cap U^+(O_K)$ is
of finite index in $U^+(O_K)$. \\

Similarly,   $\Gamma  $   intersects  $U^-(O_K)$   in   an  arithmetic
group. Hence by \cite{R 4} and \cite{V}, $\Gamma $ is arithmetic.
\end{proof}

{\bf From now on, in this section, we will assume that $K$-rank of $G$
is  ONE}. Consequently, $\fu$  has the  root space  decomposition $\fu
=\fg _{\alpha}\oplus \fg _{2\alpha}$.  Assume that $\fg _{2\alpha}\neq
0$.  Denote by $U_{2\alpha}$ the  subgroup of $G$ whose Lie algebra is
$\fg _{2\alpha}$. This is an algebraic subgroup defined over $K$. \\

It is easy to see that  the group $G_0$ whose Lie algebra is generated
by   $\fg    _{\pm   2\alpha}$   is    necessarily   semi-simple   and
$K$-simple. Moreover,  it is immediate that $S\subset  G_0$.  Note the
Bruhat decomposition  of $G$: $G=P\cup  UwP$ where $w\in N(S)$  is the
Weyl group element such that conjugation by $w$ takes $U^+$ into $U^-$
and $U_{2\alpha}$ into $U_{-2\alpha}$.  It is clear that $UwP=UwMU$ is
a Zariski open subset of $G$.

\begin{proposition} \label{highestroot} 
Suppose that  $K$ has infinitely  many units, and that  $K$-rank ($G$)
$=1$. Suppose that $\Gamma \subset G(O_K)$ is a Zariski dense subgroup
such  that $\Gamma\supset U_{2\alpha}(rO_K)$  for some  integer $r>0$.
Suppose          that         $rank$-$(G_{\infty})=\sum         _{v\in
S_{\infty}}K_v$-$rank(G)\geq 2$. Then, $\Gamma $ is of finite index in
$G(O_K)$.
\end{proposition}

\begin{proof} Let
$g=uwmv$  be an  element in  $\Gamma  \cap UwP$.   We obtain,  $\Gamma
\supset   <^g(U_{2\alpha}(rO_K)),~U_{2\alpha}(rO_K)>$.    The   Bruhat
decomposition for $g$ and  the fact that $u$ centralises $U_{2\alpha}$
shows      that      $\Gamma      \supset      ^u<U_{-2\alpha}(r'O_K),
U_{2\alpha}(r'O_K)>$ for some integer $r'$. The group $G_0$ is also of
higher real  rank, since  $S\subset G_0$ and  $K$ has  infinitely many
units. Therefore  by \cite{V}, the group generated  by $U_{\pm 2\alpha
}(r'O_K)$ is of finite index in $G_0(O_K)$ and in particular, contains
$S(r''O_K)$ for some $r''>0$. \\

We have thus seen that $\Gamma \supset ^u(S(r''O_K))$ for some integer
$r''$. Since $K$-rank of $G$ is one, the weights of $S$ acting on $\fu
$ are  $\alpha $ and  $2\alpha$. Since $S(r''O_K)$ is  infinite, there
are elements in $S(r''O_K)$ none  of whose eigenvalues acting on $\fu$
is one. Therefore, Proposition  \ref{technical} implies that $\Gamma $
is arithmetic.
\end{proof}

We  continue with  the notation  of this  subsection. There  exists an
integer $N$  such that  the units  $\theta $ of  the number  field $K$
which are  congruent to  $1$ modulo $N$,  form a  torsion-free abelian
group. Let $F$ be the  field generated by these elements. There exists
an element  $\theta \in O_K^*$ such  that for all  integers $r>0$, the
field $\Q[\theta ^r]=F$  (see Lemma \ref{exist'}).  Moreover, $S(O_F)$
is of finite index  in $S(O_K)$. We also have, 1) $F=K$  if $K$ is not
CM.  2) $F$ is totally real, $K$ totally imaginary quadratic extension
of $F$ otherwise. \\

Given  an element  $u_+\in U_{2\alpha}(O_K)\setminus  \{1\}$, consider
the subgroup $V^+$ generated  by the conjugates $^{\theta ^j}(u_+)$ of
$u_+$, as  $j$ varies over  all integers. By Lemma  \ref{exist}, there
exists      an     integer      $r$     such      that     $V^+\supset
u_+^{rO_F}\stackrel{def}{=} Exp(rO_F  log (u_+))$.  Here  $Exp$ is the
exponential map from  $Lie (U^+)$ onto $U^+$ and  $log$ is its inverse
map.   By the  Jacobson-Morozov Theorem,  there exists  a homomorphism
$f:SL(2)\ra G$ defined over $K$  such that $f\begin{pmatrix} 1 & 1\\ 0
& 1\end{pmatrix}=u_+$.  The Bruhat  decomposition shows that the image
of  the group of  upper triangular  matrices lies  in $P$.   Since all
maximal  $K$-split tori in  $P$ are  conjugate to  $S$ by  elements of
$P(K)$,  it  follows  that  there  exists  a  $p\in  P(K)$  such  that
$pf(D)p^{-1}=S$, where $D$ is the group of diagonals in $SL(2)$. Write
$p=um$ with  $u\in U$ and $m\in  M$. Now, $M$ centralises  $S$ and $u$
centralises  $u_+$ (since  $u_+$ lies  in  $U_{2\alpha}$).  Therefore,
after replacing $f$ by the map $f':x\mapsto u(f(x)u^{-1}$, we see that
$f'(D)=S$ and  $f'\begin{pmatrix} 1 & 1\\ 0  & 1\end{pmatrix}=u_+$. We
denote $f'$ by $f$ again, to avoid too much notation. \\

\begin{proposition}\label{2alpha} Suppose that $K$ has infinitely many 
units, and $G$ an absolutely  almost simple group over $K$ of $K$-rank
one. Suppose that $\fg _{2\alpha}$  is one dimensional over $K$.  Then
every arithmetic subgroup of $G(O_K)$ is virtually three generated.
\end{proposition}

\begin{proof}Let $u_+\in U_{2\alpha}(O_K)$ and $\theta \in S(O_K)$ be
as above.   Suppose that  $\gamma \in G(O_K)$  is in  general position
with respect  to $u_+$. Then, for  every $r\geq 1$,  the group $\Gamma
=<u_+^r,\theta ^r,\gamma ^r>$ is Zariski  dense. It is enough to prove
that $\Gamma$ is arithmetic.\\

By replacing $r$ by a bigger  integer if necessary, and using the fact
that  $\Z[\theta  ^r]$  has  finite   index  in  $O_F$,  we  see  that
$f\begin{pmatrix}  1  &  rO_F\\ 0  &  1\end{pmatrix}=u_+^{rO_F}\subset
\Gamma $. Write $w$ for the image of $f\begin{pmatrix} 0 & 1\\ -1 &0
\end{pmatrix}=f(w_0)$. Then, $w$ takes $U^+$ into $U^-$ under
conjugation.   Write $u_-$ for  $wu_+w^{-1}$.  Now,  $M(K)$ normalises
$U_{2\alpha}$   and  the  latter   is  one   dimensional.   Therefore,
$^m(u_+^{rO_F})=f\begin{pmatrix}    1    &     0\\    \xi    rO_F    &
1\end{pmatrix}\stackrel{def}{=}u_-^{\xi rO_F}$,  for some element $\xi
$ of  the larger field  $K$.  If $\xi  \notin F$, then  by Proposition
\ref{CM'}, the  group generated  by $u_+^{rO_F}$ and  $u_-^{\xi rO_F}$
contains   $f(\Delta)$  for   some   subgroup  of   finite  index   in
$SL(2,O_K)$. \\

Pick an element  $g\in \Gamma$ of the form  $g=uwmv$ with $u,v\in U^+$
and $m\in M$. Then,
\[\Gamma \supset <^g(u_+^{rO_F}),u_+^{rO_F}>\supset ^u
<^m(u_-^{rO_F}),  u_+^{rO_F}>=^u <u_-^{\xi rO_F},u_+^{rO_F}>\]  If $g$
is ``generic'', then $m$ is  sufficiently generic, that $\xi \notin F$
(otherwise, $^m(u_-)$ is always  $F$-rational, which lies in a smaller
algebraic group, and genericity implies  that this is not possible for
all  $m$). Then,  by the  conclusion  of the  last paragraph,  $\Gamma
\supset ^uf(S(rO_F))$ for some integer $r$.  Since $u$'s run through a
Zariski dense subset of  $U$, Proposition \ref{technical} implies that
$\Gamma $ is arithmetic.
\end{proof}

\section{Some General Results}

\subsection{Notation} Suppose $G$ is a semi-simple linear algebraic
group  which is  absolutely almost  simple and  defined over  a number
field     $K$,    with     $K$-rank($G$)$\geq    1$     and    $rank$-
$(G_{\infty})\stackrel{def}{=}\sum _{v\in S_\infty}K_v$-rank $(G) \geq
2$  (the  last  condition  says  that  $G(O_K)$  is  a  ``higher  rank
lattice'').  Let $P\subset G$  be a proper parabolic $K$-subgroup, $U$
its unipotent radical, $S\subset P$  a maximal $K$-split torus in $G$,
and $\Phi  ^+(S,P)$ the roots  of $S$ occurring  in the Lie  algebra $
\fu$ of  $U$.  Let  $\Phi ^-$ be  the negative  of the roots  in $\Phi
^+(S,P)$, and $\fu ^-=\oplus \fg _{-\alpha}$ be the sum of root spaces
with $\alpha \in \Phi ^+(S,P)$.  Then,  $\fu ^-$ is the Lie algebra of
a  unipotent  algebraic  group  $U^-$  defined over  $K$,  called  the
``opposite''  of  $U$.   Write  the  Levi  decomposition  $P=MU$  with
$S\subset M$. \\

In the  following, we will, by  restricting scalars to  $\Q$, think of
all  these  groups  $G$,  $M$,  $U^{\pm}$  as  algebraic  groups  over
$\Q$. Thus, for example, when  we say that $U^+(O_K)$ is Zariski dense
in $U^+$ we mean that $U^+(O_K)$ is Zariski dense in the complex group
$U^+(K\otimes  \C)=(R_{K/\Q}(U^+))(\C)$. With  this  understanding, we
prove  the  following  slight   strengthening  of  the  Borel  density
theorem. \\

\begin{lemma} \label{boreldense} Let $G$ be a connected 
semi-simple   $K$-simple  algebraic  group,   and  that   $G(O_K)$  is
infinite. Then, the arithmetic group  $G(O_K)$ is Zariski dense in the
complex semi-simple group $G(K\otimes \C)$.
\end{lemma}
\begin{proof} By restriction of scalars, we may assume that
$K=\Q$. Suppose that $H$ is the connected component of identity of the
Zariski closure  of $G\Z)$ in $G(\C)$. Then,  as $G(\Q)$ commensurates
$G(\Z)$, it follows that $G(\Q)$ normalises $H$.The density of $G(\Q)$
in  $G(\R)$ (weak  approximation) shows  that $G(\R)$  normalises $H$.
Clearly, $G(\R)$ is Zariski dense in $G(\C)$; hence $G(\C)$ normalises
$H$. The definition of $H$ shows  that $H$ is defined over $\Q$.  Now,
the  $\Q$-simplicity of  $G$ implies  (since $G(\Z)$  is  infinite and
hence $H$ is non-trivial) that $H=G$.
\end{proof}

The following is repeatedly used in the sequel.

\begin{lemma} \label{unipotent} Let $U$ be a unipotent group over a
number field $K$. Then, $U(O_K)$ is Zariski dense in $U(K\otimes \C)$;
moreover, if  $\Delta \subset U(O_K)$  is a subgroup which  is Zariski
dense  in $U(K\otimes  \C)$  then, $\Delta  $  is of  finite index  in
$U(O_K)$.
\end{lemma}
\begin{proof} The proof is essentially given in Theorem (2.1) of
\cite{R 5},  provided $K=\Q$. But,  by restriction of scalars,  we may
assume that $K=\Q$.
\end{proof}

Denote by  $M_0$ the  connected component of  identity of  the Zariski
closure  of $M(O_K)$ in  $M$, and  let $T_0\subset  M_0$ be  a maximal
torus  defined over  $K$.  The  groups $M$,  $M_0$ and  $T_0$  are all
defined  over  $\Q$  and  act  on  the  $\Q$-Lie  algebra  $\fu  $  of
$R_{K/\Q}U^+$  by  inner conjugation  in  $G$.   Write the  eigenspace
decomposition $\fu \otimes  \C=\oplus_{\chi \in X^*(T_0)} \fu _{\chi}$
for the action of $T_0$ on the complex lie algebra $\fu \otimes \C$.

\begin{proposition}\label{multone} 
Suppose  that each  of the  spaces $\fu  _{\chi}$ is  one dimensional.
Then  every  arithmetic  subgroup   of  $G(O_K)$  is  virtually  three
generated.
\end{proposition}

\begin{proof} Let $\mathcal U$ be the set of pairs $(m,v)\in M_0\times
\fu$  such that the  span $\sum  _{k\in \Z}\C  (^{m^k}(v))$ is  all of
$\fu$.   Then, $\mathcal  U$ is  a Zariski  open subset  of $M_0\times
\fu$. For, the  condition says that if $dim \fu  =l$, then there exist
integers $k_1,k_2,\cdots,k_l$ such that the wedge product
\[^{m^{k_1}}(v)\wedge \cdots \wedge ^{m^{k_l}}(v)\neq 0,\] 
which is a Zariski open condition.\\

Let  $\Gamma _0$  be  an  arithmetic subgroup  of  $G(K)$.  Then,  the
intersection $\Gamma _0\cap U$  is Zariski dense in $R_{K/\Q}U$.  Now,
the map $log: U\ra \fu$ is an isomorphism of varieties over $\Q$.

By assumption on $M_0$, the group $\Gamma _0\cap M_0$ is Zariski dense
in  $M_0$ (the Zariski  closure of  $\Gamma _0\cap  M_0$ is  of finite
index in $M_0$  since $\Gamma _0$ is of finite  index in $M(O_K)$, and
$M_0$ is connected).
  
By the  foregoing, we thus get  elements $m\in \Gamma  _0\cap M_0$ and
$u\in U\cap \Gamma _0$ such that $(m,log u)\in \mathcal U$. This means
that the $\Z$-span of $^{m^k}(log  u)$ as $k$ varies, is Zariski dense
in  $\fu$.   Therefore, the  group  $U_1$  generated  by the  elements
$^{m^k}(u)$  with $k\in  \Z$ is  a  Zariski dense  subgroup of  $U\cap
\Gamma _0$. Hence, by Lemma  \ref{unipotent}, $U_1$ is of finite index
in $U\cap \Gamma _0$. \\

Similarly, we can find an element $u^-\in U^-\cap \Gamma _0$ such that
the group $U_1^-$ generated by the conjugates $^{m^k}(u^-)$ with $k\in
\Z$, is of finite index in $U^-\cap \Gamma _0$. Set
\[\Gamma =<m, u,u^->\subset \Gamma _0.\]
Now, $\Gamma $ contains $<U_1,U_1^->$. By \cite{R 4} and \cite{V}, the
latter group  is of finite index  in $\Gamma _0$, hence  so is $\Gamma
$. But $\Gamma $ is three generated by construction.
\end{proof}

\begin{remark} The criterion of Proposition \ref{multone} depends on
the group $M_0$  (which is the connected component  of identity of the
Zariski closure  of $M(O_K)$)  and hence on  the $K$-structure  of the
group.   But, this  dependence is  a  rather mild  one.  However,  the
verification  that  the conditions  of  Proposition \ref{multone}  are
satisfied  is  somewhat complicated,  and  is  done  in the  next  few
sections by using the  Tits classification of absolutely simple groups
over  number  fields.  The  criterion  works  directly, when  $K$-rank
($G$)$\geq 3$,  somewhat surprisingly, for groups  of exceptional type
thanks to the analysis of  the representations of $M$ occurring in the
Lie algebra $\fu$  carried out by Langlands and  Shahidi (see \cite{L}
and \cite{Sh}). \\

However,  there  are  some  classical  groups  of  $K$-rank  $\leq  2$
(notably,  if $G$  is  of classical  type  A, C  or D  but  is not  of
Chevalley  type over  the number  field), for  which the  criterion of
Proposition \ref{multone} fails. To handle these cases, we prove below
some more lemmata of a general nature.
\end{remark}

\subsection{Notation} Let $F$ be a field of characteristic zero, and
$G$ an absolutely simple algebraic  group over $F$. Let $x\in G(F)$ be
an element of infinite order. Fix a maximal torus $T\subset G$ defined
over $F$ and $Phi$ the roots of $T$ occurring in the Lie algebra $\fg$
of $G$. We have the  root space decomposition $\fg =\ft\oplus _{\alpha
\in \Phi}\fg _{\alpha}$ with $\ft $ the Lie algebra of $T$. Now $T(F)$
is  Zariski  dense in  $T$,  hence there  exists  a  Zariski open  set
$\mathcal V \subset  T$ such that for all  $v\in T(F)\cap \mathcal V$,
the  values $\alpha  (v) ~(\alpha  \in  \Phi)$ are  all different  and
distinct from $1$. Fix $y\in T(F)\cap \mathcal V$.

\begin{lemma}\label{zariskidense} 
There is  a Zariski open set $\mathcal  U$ of $G$ such  that the group
generated by $x$ and $gyg^{-1}$ is  Zariski dense in $G$ for all $g\in
\mathcal U$.
\end{lemma}

\begin{proof} Let $H$ be a proper connected Zariski closed subgroup of 
$G$  containing (or  normalised by)  the element  $y$.  Then,  the Lie
algebra $\fh$ splits into eigenspaces for the action of $y$. Since the
values  $\alpha (y)$  are all  different (and  distinct from  $0$), it
follows   that  $\fh  =\ft\cap   \fh\oplus  \fg   _{\alpha}\cap  \fh$.
Moreover,  $\fg _{\alpha}=\fh  \cap \fg  _{\alpha}$ if  the  latter is
non-zero.  Therefore,  there exists  a proper connected  subgroup $H'$
containing $T$ which also contains  $H$ (e.g. the one with Lie algebra
$\ft  \oplus \fg  _{\alpha}$  such that  $\fh  \cap \fg_{\alpha}  \neq
0$). \\

The collection  of connected subgroups  of $G$ containing  the maximal
torus $T$  is finite since they  correspond to certain  subsets of the
finite set of  roots $\Phi$. Let $H_1,\cdots, H_n$ be  the set of {\bf
proper} connected subgroups of $G$ containing $T$. By replacing $x$ by
a power of it, we may assume that the Zariski closure $Z$ of the group
generated  by  $x$ is  connected.   Since  $G$  is simple,  the  group
generated  by $<gZg^{-1}: g\in  G>$ is  all of  $G$.  Hence,  for each
$\mu$, the set $Z_{\mu}=\{g\in G: gZg^{-1}\subset H_{\mu}\}$ is a {\bf
proper} Zariski  closed set, whence its complement  $U_{\mu}$ is open.
Therefore,  $\mathcal U=\cap  _{1\leq  \mu \leq  n}  U_{\mu}$ is  also
Zariski open. \\

Let  $g\in \mathcal  U$  and $H$  be  the connected  component of  the
Zariski  closure of  the group  generated by  $x$ and  $gyg^{-1}$.  If
$H\neq G$,  then by the first  paragraph of the proof,  there exists a
proper connected subgroup  $H'$ containing $H$ and the  torus $T$.  By
the foregoing  paragraph, $H'$  must be one  of the  $H_{\mu}$ whence,
$g\notin    U_{\mu}$,    and    so    $g\notin    \mathcal    U$,    a
contradiction. Therefore, $H=G$ and the lemma is proved.
\end{proof}

\subsection{Notation} Suppose that $G$ is an absolutely simple
algebraic group over a number  field $K$, with $K$-rank ($G$)$\geq 2$.
Let  $S$  be  a  maximal  split  torus,  and  $\Phi  (G,S)$  the  root
system. Let $\Phi ^+$ be a  system of positive roots, $\fg $ the L:lie
algebra of $G$. Let $U_0$ be  the subgroup of $G$ whose Lie algebra is
$\oplus _{\alpha >0}\fg _{\alpha}$,  and $P_0$ the normaliser of $U_0$
in $G$; then  $P_0=Z(S)U_0$ where $Z(S)$ is the  centraliser of $S$ in
$G$. Moreover, $P_0$ is a minimal parabolic $K$-subgroup of $G$. \\

Let  $\alpha  \in \Phi  ^+$  be  the highest  root  and  $\beta >0$  a
``second-highest''  root. Then,  $\gamma  =\alpha \beta$  is a  simple
root. Let $U_\alpha$ and $U_\beta$ be the root groups corresponding to
$\alpha $ and $\beta$.
\begin{proposition}\label{secondhighest} 
Let $\Gamma \subset G(O_K)$ be a Zariski dense subgroup.  Suppose that
there   exists   an   integer   $r>0$  such   that   $\Gamma   \supset
U_{\alpha}(rO_K)$ and $\Gamma \supset U_{\beta }(rO_K)$. Then, $\Gamma
$ has finite index in $G(O_K)$.
\end{proposition}

\begin{proof} This is proved in \cite{V2}. We sketch the proof, since
we  will use  this repeatedly  in the  examples.  If  $w$  denotes the
longest  Weyl group  element, then  the  double coset  $P_0wU_0$ is  a
Zariski open subset  of $G$. Hence its intersection  with $\Gamma $ is
Zariski  dense in  $G$.  Fix  an element  $g_0=p_0wu_0\in  \Gamma \cap
P_0wU_0$   and  consider  an   arbitrary  element   $g=pwu\in  PwU\cap
\Gamma$. \\

The subgroup $V=U_{\alpha }U_{\beta}$ is  normalised by all of $P$. By
assumption,  there exists  an integer  $r$ such  that  $\Gamma \supset
V(rO_K)$.    Hence   $\Gamma$   contains  the   group   $<^g(V(rO_K)),
V(rO_K)>$. By the  Bruhat decomposition of $g$, and  the fact that $V$
is   normalised  by   $P$,  we   find  an   integer  $r'$   such  that
$^g(V(rO_K))\supset                ^p(V^-(r'O_K))$               where
$V^-=U_{-\alpha}U_{w(\beta)}$  is the  conjugate of  $V$ by  $w$. Note
that $-w(\beta)$ is again a second highest root. Write $\gamma =\alpha
+w(\beta)$.  Then  $\gamma $  is a simple  root.  Moreover, it  can be
proved that the commutator subgroup $[U_{\alpha},U_{w(\beta)}]$ is not
trivial and  is all of  $U_\gamma$.  Therefore, we  get $^p(U_{\gamma}
(rO_K))\subset \Gamma $ for a  Zariski dense set of $p's$ ($r$ depends
on the  element $p$).  It  can be proved  that the group  generated by
$^p(U_{\gamma})$ is all of the  unipotent radical $U_1$ of the maximal
parabolic   subgroup  corresponding   to  the   simple   root  $\gamma
$. Consequently, for some integer $r$, $U_1(rO_K)\subset \Gamma $, and
by \cite{V2}, $\Gamma $ is arithmetic.
\end{proof}

We  will  now  deduce  a  corollary to  Lemma  \ref{zariskidense}  and
Proposition \ref{secondhighest}.

\begin{corollary} Under the notation and assumptions of this
subsection,  suppose that  every  arithmetic subgroup  $\Gamma _0$  of
$G(O_K)$ contains a 2  generated subgroup $<a,b>$ which {\bf contains}
a group of the form
\[U_\alpha (rO_K)U_\beta (rO_K)\] 
for some second highest  root $\beta$. Then, every arithmetic subgroup
of $G(O_K)$ is virtually three-generated.
\end{corollary}

\begin{proof} Given $a,b\in \Gamma _0$ such that $<a,b>\supset
(U_\alpha   U_\beta)(rO_K)$,  Lemma  \ref{zariskidense}   implies  the
existence of  an element  $c\in \Gamma $  such that the  group $\Gamma
=<a,b,c>$ generated  by $a,b,c$ is  Zariski dense in $G$  ($\Gamma _0$
itself is Zariski dense in $G$ by the Borel density theorem). Then, by
Lemma  \ref{secondhighest}, $\Gamma $  is of  finite index  in $\Gamma
_0$. \end{proof}

\newpage
\section{Groups of $K$-rank $\geq 2$}

In this section,  we verify that all arithmetic  groups of $K$-rank at
least two  are virtually three-generated.  The proof  proceeds case by
case, using the Tits classification  of algebraic groups over a number
field $K$.  In  most cases, we check that  the hypotheses of criterion
of Proposition \ref{multone}  are satisfied. In the sequel,  $G$ is an
absolutely almost simple  group defined over a number  field $K$, with
$K$-rank ($G$) $\geq 2$. The degree of $K/\Q$ is denoted $k$. \\

The  classical   groups  over  $\C$  come  equipped   with  a  natural
(irreducible) representation, which we  refer to as the {\bf standard}
representation, and denote it $St$. \\

\subsection{Groups of Inner Type A} 
In this subsection,  we consider all groups which  are inner twists of
$SL(n)$ over  $K$. By \cite{T2}, the  only such groups  are SL(n) over
number  fields or  SL(m) over  central division  algebras  over number
fields.

\subsubsection{SL(n) over number fields} $G$ is SL(n) over the number 
field $K$.  The rank assumption means  that $n\geq 3$.  Take $P$ to be
the parabolic  subgroup of  SL(n) consisting of  matrices of  the form
$\begin{pmatrix} g &  x \\ 0 & det  g^{-1} \end{pmatrix}$, where $g\in
GL_{n-1}$,       $x=\begin{pmatrix}      x_1\\       x_2\\      \cdots
\\x_{n-1}\end{pmatrix}$ is a  column vector of size $n-1$,  and $0$ is
the $1\times (n-1)$  matrix whose entries are all  zero. The Levi part
$M$ of $P$ may be taken to  be $GL(n-1)=\{\begin{pmatrix} g & 0 \\ 0 &
det(g^{-1})\end{pmatrix}:g\in GL_{n-1}\}$.   Recall that $M_0$  is the
connected   component  of   identity   of  the   Zariski  closure   of
$M(O_K)$. Hence $M_0$ contains the subgroup $H=SL(n-1)$. Take $T_0$ to
be the  diagonals in  $SL(n-1)$. The unipotent  radical of $P$  is the
group $\begin{pmatrix} 1 & x\\ 0  & 1 \end{pmatrix}$ with $x$ a column
vector as before. As  a representation of $GL(n-1)(\C)$ (and therefore
of  $H(\C)=SL(n-1)(\C)$),  the  Lie   algebra  $\fu$  is  nothing  but
$st\otimes  det   $,  the  standard  representation   twisted  by  the
determinant. Restricted to the torus $T_0$, thus $\fu$ is the standard
representation,  and  is  hence  multiplicity  free.   By  Proposition
\ref{multone}, it  follows that every arithmetic  subgroup of $G(O_K)$
is virtually three-generated.

\subsubsection{SL(m) over division algebras} $G =SL_m(D)$ where $D$ is
a central division algebra over the number field $K$, of degree $d\geq
2$. The  rank assumption  means that $m\geq  3$. Consider  the central
simple algebra  $D\otimes _{\Q}\R$,  denoted $D\otimes \R$  for short.
Then, $D\otimes \R$ is a product of copies of $M_d(\C)$, $M_{d/2}({\bf
H})$,  and  $M_d(\R)$ where  ${\bf  H}$  is  the division  algebra  of
Hamiltonian quaternions. We consider four cases.\\

{\it Case 1. $D\otimes \R\neq {\bf H}\times \cdots {\bf H}$}.
  
Then, $SL_1(D\otimes \R)$ is  a non-compact semi-simple Lie group with
either $SL_d(\C)$  or $SL_d(\R)$ or $SL_{d/2}({\bf H})$  (the last can
happen only if  $d\geq 3$ is even) as a  non-compact factor. Then, the
Zariski closure  of the  arithmetic subgroup $SL_1(O_D)$  of $SL_1(D)$
(for some order  $O_D$ of $D$) is the  noncompact group $SL_1(D\otimes
\R)$.\\

Take  $P$ to  be the  parabolic subgroup  (with the  obvious notation)
$P=\begin{pmatrix} GL_1(D) & *  \\ 0 & GL_{m-1}(D)\end{pmatrix}$, with
unipotent  radical $U=\begin{pmatrix}  1  & M_{1\times  m}(D)  \\ 0  &
1\end{pmatrix}$ where $M_{1times m}(D)$ denotes the spaces of $1\times
m$ matrices with entries in  the division algebra $D$. The group $M_0$
obviously contains  (from the observation  in the last  paragraph) the
group $H=\begin{pmatrix}SL_1(D) &  0 \\ 0 & SL_{m-1}(D)\end{pmatrix}$,
with $H(K\otimes \C)=[SL_d(\C)\times SL_{d(m-1)}(\C)]^k$. Let $T_0$ be
the   product  of   the  diagonals   in  each   copy   of  $SL_d\times
SL_{d(m-1)}$. As a representation of $H$, the Lie algebra $\fu$ of $U$
is the direct sum $\oplus  \C^d\otimes \C ^{(m-1)d})^*$, where the sum
is over each copy of $SL_d\times SL_{(m-1)d}$.  $\C^d$ is the standard
representation of  $SL_d$ and $*$ denotes  its dual. It  is then clear
that as a representation of the product diagonal torus $T_0$, $\fu$ is
multiplicity  free. Hence the  criterion of  Proposition \ref{multone}
applies.   Every arithmetic  subgroup of  $G(O_K)$ is  virtually three
generated.\\

{\it Case  2.  $D\otimes \R={\bf  H}\times \cdots \times {\bf  H}$ but
$m\geq 4$.}

This can happen  only if $d=2$, and $D\otimes  \R={\bf H}^k$. Take $P$
to be the parabolic subgroup  $P=\begin{pmatrix}SL_{m-2}(D) & * \\ 0 &
SL_2(D)\end{pmatrix}$ and  denote its  unipotent radical by  $U$, with
$U=\begin{pmatrix}1  & M_{(m-2)\times  2}(D) \\  0  & 1\end{pmatrix}$.
Then, as before, $M_0$ contains $H=\begin{pmatrix}SL_{m-2}(D) & 0 \\ 0
& SL_2(D)\end{pmatrix}$.  Then, $H(K\otimes \C)=[SL_{(m-2)2}(\C)\times
SL_4(\C)]^k$.  As a representation of  $H(K\otimes \C)$, $\fu $ is the
$k$-fold direct sum $\C^{(m-2)2}\otimes  (\C ^4)^*$.  Let $T_0$ be the
product of the diagonals in  $H(K\otimes \C)$.  Then, it is clear that
$\fu$  is  multiplicity  free   as  a  representation  of  $T_0$.   By
Proposition  \ref{multone} every  arithmetic subgroup  of  $G(O_K)$ is
virtually three-generated. \\

{\it Case 3.  $D\otimes \R={\bf H}\times \cdots \times {\bf H}$, $m=3$
  but $k\geq 2$}.

In  this  case  it  turns   out  that  the  criterion  of  Proposition
\ref{multone} fails ( we will not  prove that it fails), so we give an
ad  hoc  argument  that  every  arithmetic subgroup  of  $SL_3(D)$  is
virtually three-generated. \\

Since $D\otimes  \R={\bf H}^k$, it  follows that $K$ is  totally real.
Since   $k\geq  2$,  $K$   has  infinitely   many  units.    By  lemma
\ref{finite}, for every subgroup $\Delta $ of finite index in $O_K^*$,
$\Q[\Delta]=K$. By Lemma \ref{exist},  there exists an element $\theta
\in \Delta $  such that $\Q[\theta ^r]=K$ for  all $r\geq 1$. Consider
the $3\times 3$ - matrix $m=  \begin{pmatrix} \theta ^{k_1} & 0 & 0 \\
0 &  \theta ^{k_2}  & 0  \\ 0 &  0 &  \theta ^{-k_1-k_2}\end{pmatrix}$
which lies in $SL_3(O_D)$ for some order $O_D$ in $D$. \\

Consider the following matrices in $SL_3(D)$ given by $u=
\begin{pmatrix} 1 & 1 & 1 \\ 0 & 1 & x \\ 0 & 0 & 1\end{pmatrix}$, and
$u^-=  \begin{pmatrix}  1  &  0 &  0  \\  y  &  1  &  0  \\ z  &  t  &
1\end{pmatrix}$,  where   $x,y,z,t$  are  elements   of  the  division
algebra,no two  of which commute. We  may assume that they  lie in the
order $O_D$.   We will prove that  for every $r>0$,  the group $\Gamma
=<m^r,u^r,(u^-)^r>$  generated by  the $r$-th  powers of  $m,u,u^-$ is
arithmetic. This will prove  that every arithmetic subgroup of $GO_K)$
is  virtually three  generated.  We  use the  following  notation.  If
$i,j\leq 3$  , $i\neq  j$ and $w$  is an  element of $O_D$,  denote by
$x_{ij}^{O_Kw}$  the  subgroup $1+cwE_{ij}$  where,  $c$ runs  through
elements of $O_K$;  $E_{ij}$ is the matrix whose  $ij$-th entry is $1$
and  all other  entries are  zero.  We  also  write $x_{ij}^{O_Kw}\leq
\Gamma  $   to  say   that  for  some   integer  $r'$,   the  subgroup
$x_{ij}^{r'O_Kw}$ is contained in $\Gamma $. \\

For ease of  notation, we replace the r-th powers  of $m,u,u^-$ by the
same  letters $m,u,u^-$; this  should cause  no confusion.   The group
$<m^l(u): l\in \Z>$ virtually contains the subgroups (by the choice of
$\theta$; see  Lemma \ref{exist}) $x_{12}^{O_K}$  , $x_{13}^{O_K}$ and
$x_{23}^{xO_K}$.   Similarly, $<m^l(u^-):l\in \Z>$  virtually contains
$x_{21}^{O_Ky},  x_{31}^{O_Kz}$ and  $x_{32}^{O_Kt}$. Since  $\Gamma $
contains   all   these  groups,   by   taking   commutators,  we   get
$x_{12}^{O_Kt}=[x_{13}^{O_K},x_{32}^{tO_K}]\leq  \Gamma $.  Similarly,
$x_{13}^{O_Kx}=[x_{12}^{O_K},x_{23}^{O_Kx}]\leq      \Gamma,$      and
$x_{13}^{O_Ktx}=[x_{12}^{O_Kt},x_{23}^{txO_K}]\leq   \Gamma   $.    By
taking suitable commutators, we obtain
\[x_{12}^{O_K+O_Kt+O_Kx+O_Ky}=\leq \Gamma. \] Since up to subgroups of
finite  index $O_D=OK+O_Kt+O_Kx+O_Ky$,  we  see that  $x_{13}(O_D)\leq
\Gamma  $,  and similarly,  $x_{ij}(O_D)$  for  all  $ij$ with  $i\neq
j$. Therefore, $\Gamma \supset U(O_K)$ and $U^-(O_K)$ for two opposing
maximal unipotent subgroups  of $G$. By \cite{R 4},  $\Gamma $ is then
an arithmetic group. \\

{\it Case 4. $D\otimes \R={\bf H}^k$, $m=3$ and $k=1$.}

The assumptions mean that  $K=\Q$, and $D\otimes \R={\bf H}$. Consider
the elements $m_0= \begin{pmatrix} a  & b & 0 \\ c & d &  0 \\ 0 & 0 &
1\end{pmatrix}$, $u_0= \begin{pmatrix} 1 & 0 & 1 \\ 0 & 1 & 0 \\ 0 & 0
& 1\end{pmatrix}=x_{13}$, $u_0^-= \begin{pmatrix} 1 & 0 & 0 \\ 0 & 1 &
0  \\ 0  &  1 &  1\end{pmatrix}=x_{32}$.   We assume  that the  matrix
$\begin{pmatrix}  a &  b \\  c &  d\end{pmatrix}$ is  ``generic''.  In
particular,  assume that  $c\notin  \Q$ and  that  $e=ca+dc$ does  not
commute    with    $c$.    Fix    $r\geq    1$    and   put    $\Gamma
=<m_0^r,u_0^r,(u_0^-)^r>$. By  arguments similar to the  last case, it
is  enough  to prove  that  $\Gamma $  is  an  arithmetic subgroup  of
$SL_3(O_D)$ for some order $O_D$ of the division algebra $D$. \\

We have  $x_{13}^{\Z}\leq \Gamma $ and $x_{32}^{\Z}\leq  \Gamma $.  By
taking commutators, we get $x_{12}^{\Z}\leq \Gamma $. \\

The conjugate $^{m^0}(u_0)=\begin{pmatrix} 1 & 0 & a \\ 0 & 1 & c \\ 0
& 0 & 1\end{pmatrix}$.  Hence we get $\begin{pmatrix} 1 & 0 & a \\ 0 &
1  &  c \\  0  &  0 &  1\end{pmatrix}^{\Z}\leq  \Gamma  $.  By  taking
commutators with $x_{12}^{\Z}$, we then get $x_{13}^{\Z[c]}\leq \Gamma
$.   Taking  commutators with  $x_{32}^{\Z}\leq  \Gamma  $, we  obtain
$x_{12}^{\Z[c]}\leq \Gamma $ as well. \\

Consider $\begin{pmatrix} a & b \\ c & d \end{pmatrix}
\begin{pmatrix} a \\ c \end{pmatrix}=
\begin{pmatrix} a^2+bc \\ ca+dc \end{pmatrix}= 
\begin{pmatrix} a' \\ e \end{pmatrix}$. Clearly, $^{m_0^2}(u_0)=
\begin{pmatrix} 1 & 0 & a'\\ 0 & 1 & e \\0 & 0 & 1 \end{pmatrix}.$ By
the  argument  of  the  last  paragraph,  taking  commutators  of  its
conjugate  with  $x_{12}^{\z[c]}$  we obtain  $x_{13}^{\Z[e]\Z[c]}\leq
\Gamma $. Since  $e$ and $c$ do not commute and  $D$ has dimension $4$
over $\Q$, it follows that $\Z[e]\Z[c]$ is of finite index in an order
$O_D$  of   $D$.   Therefore,  $x_{13}^{O_D}\leq   \Gamma  $.   Taking
commutators    with    $x_{32}^{\Z}\leq    \Gamma   $,    we    obtain
$x_{12}^{O_D}\leq  \Gamma $ (and  $x_{13}^{O_D}\leq \Gamma  $).  Thus,
$\Gamma $ intersects the unipotent radical (consisting of $x_{12}$ and
$x_{13}$ root groups) of a parabolic subgroup of $G$. Clearly, $\Gamma
$ is Zariski dense. Therefore, by \cite{V2}, $\Gamma $ is arithmetic.

\subsection{Groups of outer type A}

Suppose that $K$ is a number  field of degree $k\geq 1$ over $\Q$. Let
$E/K$  be a  quadratic  extension  and $\sigma  \in  Gal(E/K)$ be  the
non-trivial element.   Suppose that $D$ is a  central division algebra
over $E$ of degree $d\geq  1$ (as usual $d^2=dimension (D/E)$). Assume
there is  an involution $*$  on $D$ such  that its restriction  to the
centre $E$  coincides with $\sigma$. If  $N\geq 1$ is  an integer, and
$g\in M_N(D)$, with $g=(g_{ij})$ is an $N\times N$ matrix with entries
in $D$,  then define $g^*$ as  the matrix with $ij$-th  entry given by
$g^*_{ij}=(g_{ij})^*$.   Thus, $M_N(D)$  gets an  involution $g\mapsto
(^tg)^*$. \\

Fix  an integer  $m\geq  0$. Consider  the $(m+4)\times  (m+4)$-matrix
$h=\begin{pmatrix} 0_{2\times 2} &  0_{2\times m} & \begin{pmatrix}0 &
1 \\ 1 & 0 \end{pmatrix}\\ 0_{m\times 2} & h_0 & 0_{m \times 2}\\
\begin{pmatrix} 0 & 1\\ 1 & 0\end{pmatrix} & 0_{2\times m} &
0_{2\times 2}  \end{pmatrix}$, where $0_{p\times q}$  denotes the zero
matrix of the relevant size.  Here $h_0$ is a non-singular $m\times m$
matrix  with entries  in  $D$  such that  $(^t  h_0)^*=h_0$. Then  $h$
defines a Hermitian form with  respect to $*$ on the $m+4$ dimensional
vector space over $D$. The algebraic group we consider is of the form
\[G=SU_{m+4}(h,D)=\{ g\in SL_{m+4}(D): (^tg)^*hg=h \}.\]
Then $G$ is  an absolutely simple algebraic group  over $K$. Since $h$
contains  {\bf  two}  copies  of  the  ``hyperbolic''  Hermitian  form
$J=\begin{pmatrix}  0  &  1\\  1  & 0\end{pmatrix}$  it  follows  that
$K$-rank ($G$)$\geq  2$. From the classification  tables of \cite{T2},
these $G$  are the only outer forms  of type $A$ of  $K$-rank at least
two.\\

Arithmetic  subgroups $\Gamma _0$  of $G$  are commensurate  to $G\cap
GL_{m+4}(O_D)$  for   some  order   $O_D$  of  the   division  algebra
$D$. Consider the subgroup
\[H=\{\begin{pmatrix}g & 0 & 0\\ 0 & 1_m & 0\\ 0 & 0 &
  J[(^tg)^*]^{-1}J^{-1}  \end{pmatrix}:   g\in  SL_2(D)\}.   \]  Since
$H(K\otimes \R)=SL_2(D\otimes \R)$ is non-compact, it follows from the
Borel  density   theorem  (see  (.))   that   $\Gamma  _0\cap  H\simeq
SL_2(O_D)$. Moreover, $\Gamma _0\cap H$  is Zariski dense in the group
$H(K\otimes \C)=[SL_{2d}(\C)\times  SL_{2d}(\C)]^k$.  The intersection
of $G$  with diagonals is at  least two dimensional, and  is a maximal
$K$-split torus $S$,  if $h_0$ is suitably chosen  (that is, split off
all the hyperbolic forms in $h_0$ in the same way as was done for {\it
two} hyperbolic forms for $h$). \\

With respect  to $S$, the  intersection of unipotent  upper triangular
matrices with $G$ yields a maximal unipotent subgroup $U_0$ of $G$ and
the roots of $S$ occurring in the Lie algebra $\fu _0$ of $U_0$ form a
system $\Phi ^+$ of positive roots.  If $\alpha $ and $\beta $ are the
highest  and  a second  highest  root in  $\Phi  ^+$,  then the  group
$U_\alpha  (O_K)U_\beta   (O_K)$  is   contained  in  the   group  $U=
\{\begin{pmatrix} 1 &  0 & x \\ 0 &  1 & 0 \\ 0  & 0 & 1\end{pmatrix}:
x\in M_2(O_D),  (^tx)^*+JxJ=0\}$.  Now,  as a module  over $H(K\otimes
\C)=[SL_{2d}(\C)\times SL_{2d}(\C)]^k$, the Lie algebra $Lie U(\C)$ is
isomorphic  to  $[\C^{2d}\otimes   (\C^{2d})^*]^k$  and  has  distinct
eigenvalues for the diagonals $T_H$ in $H(K\otimes \C)$ (thought of as
a product of copies of $SL_{2d}(\C)$).  Therefore, by section 3, there
exist $m_0\in \Gamma _0\cap H$, $u_0\in U\cap \Gamma _0$ such that the
group  generated  by   the  conjugates  $\{^{m_0^j}(u_0):  j\in  \Z\}$
contains the group $U(rO_K)$ for  some integer $r$.  Now the criterion
of Proposition  \ref{secondhighest} says  that there exists  a $\gamma
_0\in \Gamma _0$ such that the three-generated group $\Gamma =<\gamma,
m_0, u_0>$ is of finite index in $\Gamma _0$. \\

\subsection{Groups of type B and inner type D} (i.e. type $^1D^1_{n,r}$)

In this subsection, we consider  groups of the form $G=SO(f)$ with $f$
a non-degenerate quadratic  form in $n$ variables over  $K$, $n\geq 5$
(and $n\geq 8$ if $n$ is even). Assume that $f$ is a direct sum of two
copies  of a hyperbolic  form and  another non-degenerate  form $f_2$:
$f=\begin{pmatrix} 0 & 1\\ 1 & 0\end{pmatrix}\oplus
\begin{pmatrix} 0 & 1\\ 1 & 0\end{pmatrix}\oplus
f_2$. \\

Put $f_1=\begin{pmatrix}  0 & 1\\ 1 &  0\end{pmatrix}\oplus f_2$. Then
$f=\begin{pmatrix}  0  & 1\\  1  &  0\end{pmatrix}\oplus f_1$.   Then,
$K$-rank ($G$) $\geq 2$. Consider the subgroup $P= \{\begin{pmatrix} a
&  x   &  -\frac{x^tx}{2}\\  0   &  SO(f_1)  &   -^tx  \\  0  &   0  &
a^{-1}\end{pmatrix}: a\in  {\bf G}_m, x\in K^{n-2}\}$.  Then  $P$ is a
parabolic subgroup of $G$ with unipotent radical $U= \{\begin{pmatrix}
1 & x & -\frac{x^tx}{2}\\ 0 & 1 & -^tx \\ 0 & 0 & 1\end{pmatrix}: x\in
K^{n-2}\}$. Now, the group $SO(f_1)$ is isotropic over $K$ since $f_1$
represents  a  zero.  Moreover,  since  $n-2\geq  3$,  $SO(f_1)$ is  a
semi-simple algebraic group over  $K$. Hence $SO(f_1)(O_K)$ is Zariski
dense  in $SO(f_1)(K\otimes  \C)$.  Consequently,  $M_0$  contains the
subgroup $SO(f_1)$. Moreover, as a representation of $SO(f_1)(K\otimes
\C)=SO(n-2)(\C)^k$, the  Lie algebra $\fu(K\otimes \C)$ of  $U$ is the
standard  representation $St  ^k$.  Clearly,  for a  maximal  torus in
$SO(n-2)(\C)$,    the   standard   representation    is   multiplicity
free. Therefore,  the criterion of  Proposition \ref{multone} applies:
every arithmetic subgroup of $G(O_K)$ is virtually three-generated.

\subsection{Groups of type C and the rest of the Groups of type D}

\subsubsection{$G=Sp_{2n}$ over $K$ with $n\geq 3$.} 

Denote by
\[\kappa =\begin{pmatrix} 0 & 0 & \cdots & 0 & 1\\ 0 & 0 &
\cdots & 1 & 0\\ 0 & \cdots &\cdots &\cdots & 0\\ 1 & 0 & \cdots & 0 &
0\end{pmatrix}\] the $n\times n$ matrix all of whose entries are zero,
except for  the anti-diagonal ones, which  are all equal  to one.  Let
$J=\begin{pmatrix} 0_n & \kappa  \\ -\kappa & 0_n\end{pmatrix}$ be the
non-degenerate  $2n\times  2n$  skew  symmetric  matrix.   Define  the
symplectic  group $G=Sp_{2n}=\{ g\in  SL_{2n}: ^tgJg=J\}$.   The group
$P=\{\begin{pmatrix} g & 0 \\ 0 & \kappa ^tg\kappa ^{-1}
\end{pmatrix}\begin{pmatrix} 1 & x \\0 & 1\end{pmatrix}: x+\kappa ^tx
\kappa =0, g\in  GL_n \}$ is a parabolic subgroup.   Denote by $M$ the
Levi subgroup of  $P$ such that $x=0$.  Then, it is  easy to show that
$M_0\supset H\simeq  SL_n=\{\begin{pmatrix} g  & 0 \\  0 &  \kappa ^tg
\kappa  ^{-1}\end{pmatrix}\}. $ As  a representation  of $H$,  the Lie
algebra  $\fu$ of  the unipotent  radical  $U$ of  $P$ is  seen to  be
isomorphic to $S^2(\C^n)$, the  second symmetric power of the standard
representaqtion of  $H=SL_n$. Therefore, with respect  to the diagonal
torus  $T_H$ of $H$,  the representation  $\fu$ is  multiplicity free.
Therefore, by Proposition  \ref{multone}, every arithmetic subgroup of
$Sp_{2n}(O_K)$ is three-generated.
  
\subsubsection{Other Groups of type C and D}
In this  subsection, we will consider  all groups of type  $C$ or $D$,
which are not  covered in the previous subsections. Let  $D$ be a {\bf
quaternionic} division algebra over the number field $K$. Let $\sigma$
be an involution  (of the first kind)  on $D$.  In the case  of type C
(resp.   type  D),  assume  that  the space  $D^{\sigma}$  of  $\sigma
$-invariants in $D$ is  one dimensional (resp. three dimensional) over
$K$.  Let  $m\geq 0$  be an integer.   Consider the  $m+4$ dimensional
matrix  $h=\begin{pmatrix}0_2  & 0  &  \begin{pmatrix}0  &  1 \\  1  &
0\end{pmatrix}\\ 0 & h_0 & 0 \\ \begin{pmatrix} 0 & 1 \\ 1 & 0
\end{pmatrix} & 0 & 0 \end{pmatrix}$, where $h_0$ is a non-singular
matrix with entries  in $D$, such that $^t\sigma  (h_0)=h_0$.  We will
view  $h$ as  a non-degenerate  form on  $D^{m+4}\times  D^{m+4}$ with
values  in $D$,  which is  hermitian  with respect  to the  involution
$\sigma$. The algebraic group which we consider is the special unitary
group of this hermitian form:  $G=SU(h)$ - an algebraic group over $K$
(if $D^{\sigma}$  is three dimensional, then  $G$ is of  type $^1D$ or
$^2D$ according as the discriminant  of $h$ is $1$ or otherwise). With
this  choice of  $h$,  it is  immediate  ($h$ has  two  copies of  the
hyperbolic form: cf.  the subsection  on groups of outer type A), that
$K$-rank ($G$) $\geq 2$.\\

Since We  needed $K$-rank ($G$) $\geq  2$ we split  off two hyperbolic
planes from $h$. The form $h_0$ may have more hyperbolic planes in it;
after splitting  these off  in a  manner similar to  that for  $h$, we
obtain a form $h_1$ which is anisotropic over $K$. We will assume that
$h_0$  is  of this  type.   Then, the  intersection  of  $G$ with  the
diagonals is  a maximal $K$-split torus  $S$ in $G$. The  roots of $S$
occurring in the  group of unipotent upper triangular  matrices in $G$
form a positive  system $\Phi ^+$.  Choose $\alpha  $ the highest root
and a second highest root $\beta$ in $\Phi ^+$. \\

The  group $U_\alpha  U_\beta $  is contained  in the  unipotent group
$U=\{\begin{pmatrix}1_2 & 0 & x \\ 0 & 1 & 0 \\ 0 & 0 & 1
\end{pmatrix}: x suitable\}$ which is the unipotent radical of a
parabolic subgroup. Set $H=\{\begin{pmatrix} g & 0 & 0 \\ 0 & 1 & 0 \\
0  &  0 &  J^tgJ^{-1}\end{pmatrix}  :  g\in  SL_2(D)\}$.  Then,  $M_0$
contains  $H$.   Let  $\Gamma   _0\subset  G(O_K)$  be  an  arithmetic
subgroup. Then,  there exist $m_0\in  H(O_K)\cap \Gamma _0$  and $u\in
U(O_K)\cap \Gamma _0$  such that the group generated  by $m_0$ and $u$
(denoted as usual by  $<m_0,u>$) intersects $(U_\alpha U_\beta )(O_K)$
in  a subgroup  of finite  index. By  Lemma  \ref{zariskidense}, there
exists  an   element  $\gamma  \in   \Gamma  _0$  such   that  $\Gamma
=<m_0,u,\gamma >$ is Zariski dense in $G(K\otimes \C)$. By Proposition
\ref{secondhighest},  $\Gamma  $  has  finite index  in  $\Gamma  _0$:
$\Gamma _0$ is virtually three-generated. \\
   
\subsection{The Exceptional Groups}

In  this  subsection,  we  prove  Theorem  1 for  all  groups  $G$  of
exceptional type of $K$-rank $\geq 2$. In each of these cases, we will
locate a simple $K$-root in the Tits -Dynkin diagram of $G$, such that
the Levi subgroup (actually the  group $M_0$ contained in the Levi) of
the  parabolic  group corresponding  to  the  simple  root contains  a
subgroup  $H$ with  the following  property.  $H(K\otimes  \C)$  has a
maximal  torus  $T_H$  whose  action  on the  Lie  algebra  $\fu  =Lie
(U)(K\otimes \C)$ is  multiplicity free. By Proposition \ref{multone},
this  implies that every  arithmetic subgroup  of $G(K)$  is virtually
three-generated. The notation is as in \cite{T2}. \\

\subsubsection{the groups $^3D^2_{4,2}$ and $^6D^2_{4,2}$.} 

In the  Tits diagram, there is  one simple circled  root $\alpha$, and
three other simple roots  which are circled together.  The semi-simple
part of  the Levi  is therefore $K$-simple,  and hence  contains (over
$\C$),  the  group   $SL_2(\C)^3$  (three-fold  product  of  $SL(2)$).
According to  \cite{L} and \cite{Sh}, the representation  $\fu$ is the
direct  sum of  $St\otimes St\otimes  St$ and  $1\otimes  St\otimes 1$
($St$  is the  standard representation  and $1$  is the  trivial one).
This  is multiplicity free  for the  product of  the diagonals  in the
group $SL(2)^3$.

\subsubsection{Groups of type $E_6$}

There are three groups of inner type $E_6$ with $K$-rank $\geq 2$.

{\it Case 1. $G=^1E_{6,2}^{28}$}. The extreme left root in the diagram
is  circled.   Since  its  $K$-rank  is  $\geq 1$,  the  Levi  of  the
corresponding maximal  parabolic subgroup is  non-compact. Then, $M_0$
contains $SO(10)$ over $\C$. According to \cite{L}, the representation
on $\fu$  is one of the $\frac{1}{2}$-spin  representation of $SO(10)$
and has distinct characters for the maximal torus.  \\

{\it Case  2.  $G= ^1E_{6,2}^{16}$}.   Over the number field  $K$, the
diagram is  that of  $^1E_{6,2}^{16}$. The root  in the middle  of the
diagram  is circled.   However, over  any archimedean  completion, the
diagram can only be the split form ($^1E_{6,2}^{16}$ can not transform
into $^1E_{6,2}^{28}$  over $\R$ or  $\C$).  Consequently, $M(K\otimes
\R)$  contains $SL_3\times SL_2\times  SL_3$, whence  $M_0= SL_3\times
SL_2\times  SL_3$.   According  to  \cite{Sh}, the  representation  of
$M_0(\C)$ on $\fu$ is the  direct sum of $St _{SL_3}\otimes St _{SL_2}
\otimes \wedge ^2  St_{SL_3}$ (from now on we  will drop the subscript
$SL_3$ or  $SL_2$ for  ease of notation),  $\wedge ^2 St  \otimes Triv
\otimes  St$ and  $Triv \otimes  St \otimes  Triv$. It  is  clear that
restricted to  the product of the diagonals  in $SL_3\times SL_2\times
SL_3$, the representation $\fu $ has multiplicity one. \\
  
{\it Case 3. $G= ^1 E_{6,6}^0$}. The same $M_0$ as in Case 2 works, to
prove multiplicity one for the torus. \\

Now  consider the groups  of outer  type $E_6$  of $K$-rank  $\geq 2$.
These   are  $^2E_{6,2}^{16   ^{'}}$,   $^2E_{6,2}^{16  ^{''}}$,   and
$^2E_{6,4}^2$. In all these, the  root at the extreme left is circled.
Then, since  $K$-rank ($M$) $\geq 1$, it  follows that $M_0(\C)\supset
SL_6$.  The representation on $\fu$ is (by \cite{L}, page 49, (x)) the
direct sum  of $Triv$  and $\wedge ^3(St)$  and the diagonal  torus in
$SL_6$ has multiplicity one for its action on $\fu$. \\

\subsubsection{Groups of type $E_7$}

There  are four  groups of  type $E_7$  over a  number field  $K$ with
$K$-rank   $\geq  2$.    They   are  $E_{7,2}^{31}$,   $E_{7,3}^{28}$,
$E_{7,4}^9$ and  $E_{7,7}^0$. In  all these, the  root on  the extreme
right  is circled,  and $M$  has  $K$-rank $\geq  1$. Hence  $M_0(\C)$
contains  the  semi-simple  part  of  the Levi  group  $M$.   This  is
$SO(12)$. According to \cite{L},  the representation $\fu$ of $SO(12)$
is $triv\oplus  \frac {1}{2} spin$ which has  distinct eigenvalues for
the torus in $SO(12)$.

\subsubsection{Groups of type $E_8$} 

The groups with $K$-rank  $\geq 2$ are $E_{8,4}^{28}$ and $E_{8,8}^0$.
Consider  the  root  on  the   extreme  right  in  the  diagram.   The
corresponding  $M$  has semi-simple  (actually  simple) part  $SO(14)$
which is  isotropic over $K$.  The representation  $\fu$, according to
\cite{L}, is  $\frac{1}{2}$-spin $\oplus St$ and  has multiplicity one
for the maximal torus of $SO(14)$.

\subsubsection{The groups $F_4$} 

There  is only  one  $K$-rank  $\geq 2$  group,namely  the split  one,
denoted  $F_{4,4}^0$.   Take  the  root  on the  extreme  left.   Then
$M_0\supset    SO(7)$.    The    representation   is    $triv   \oplus
\frac{1}{2}-spin$  and is  multiplicity  free for  the  action of  the
maximal torus in $SO(7)$.

\subsubsection{ Groups of type $G_2$}. The only group is $G_{2,2}^0$,
the split  form. For  the root  on the extreme  left, the  group $M_0$
contains $SL(2)$  and the representation $\fu$ is  $Triv \oplus Sym^3$
which has distinct eigenvalues for  the action of the maximal torus in
$SL(2)$.

This completes the proof of Theorem 1 for groups of $K$-rank $\geq 2$.

\newpage

\section{Classical Groups of Rank One} 

The case  of groups $G$  such that $K$-rank  ($G$) $= 1$  and $\R-rank
(G_\infty)\geq 2$ is much more involved. We will have to consider many
more cases, both for classical  and exceptional groups. In some cases,
we will  have to  supply ad hoc  proofs, because the  general criteria
established in the previous sections do not apply.
\subsection{Groups of inner type A} 

The assumptions imply that $G=SL_2(D)$ where $D$ is a central division
algebra over the number field $K$.\\

{\it Case 1. $D=K$}. The assumption that $\R-rank (G_\infty)\geq 2$ is
equivalent  to $r_1+r_2\geq  2$.  Therefore,  $K$ has  infinitely many
units. This has been covered in Section (2.1) on $SL(2)$.\\

{\it Case  2. $D\neq K$,  $D\otimes _{\Q}\R \neq {\bf  H}\times \cdots
\times{\bf H}$}.   Here ${\bf H}$  denotes the algebra  of Hamiltonian
quaternions.  Consider the parabolic subgroup $P=\{\begin{pmatrix} g &
0 \\  0 & h\end{pmatrix} \begin{pmatrix}1  & x \\  0 & 1\end{pmatrix}:
g,h\in  GL_1(D), Det  (gh)=1,  x\in  D\}$. Let  $U$  be its  unipotent
radical.  The assumption on $D$  means that $SL_1(D\otimes \R)$ is not
compact,  and  $SL_1(D\otimes  _{\Q}\C)$  contains  $SL_1(O_D)$  as  a
Zariski  dense  subgroup  (Proposition \ref{boreldense}).   Therefore,
$M_0$    contains    the    subgroup    $M_1$    with    $M_1(K\otimes
\C)=[SL_2(\C)\times  SL_2(\C)]^k$. As  a representation  of $M_1(\C)$,
the Lie  algebra $\fu =(Lie  U)(K\otimes \C)$ is  $[St\otimes St^*]^k$
and  is  multiplicity  free  for  the  action  of  the  maximal  torus
($2k$-fold product  of the diagonals in  $M_1(\C)$).  Therefore, every
arithmetic subgroup of $G(K)$ is virtually three-generated. \\
  
{\it Case  3. $D\neq K$ and $D\otimes  _{\Q}\R={\bf H}^k$}. Therefore,
$K$ is totally  real of degree $k$ over  $\Q$. The assumption $\R-rank
(G_\infty)\geq  2$  means  that  $K\neq  \Q$.  Let  $P$,  $M$  be  the
parabolic  subgroup and its  Levi subgroup  in Case  2 of  the present
subsection.  Fix $m=\begin{pmatrix} \alpha & 0\\ 0 & \beta
\end{pmatrix}\in M(K)$ such that $\delta =\alpha \beta ^{-1}$ does not
lie  in  $K$. Since  $D\otimes  \R={\bf  H}^k$,  it follows  that  the
extension $K(\delta )/K$ is a CM extension. Fix $u_+=\begin{pmatrix} 1
& 1\\ 0  & 1 \end{pmatrix}$ and $m_0=\begin{pmatrix} \theta  & 0\\ 0 &
\theta ^{-1} \end{pmatrix}$  where $\theta \in O_K^*$ is  chosen as in
Lemma \ref{exist}.   Fix $\gamma \in G(O_K)$ in  general position with
respect to  $u_+$ and  $m_0$. Then, for  every integer $r$,  the group
$\Gamma =<u_+^r,m_0^r,\gamma  ^r>$ generates a  Zariski dense subgroup
of $G$ (see Lemma \ref{zariskidense}).  We will show that $\Gamma $ is
arithmetic,  proving  that  every  arithmetic subgroup  of  $G(K)$  is
virtually 3-generated (since arithmetic  groups contain a group of the
form $<\gamma ^r,m_0^r,u_+^r>$ for some integer $r$).\\

Since $\Gamma  $ contains $\theta ^r$  and $u_+$, it  follows that for
some integer $r'$, $\Gamma $ contains the group $V^+=\begin{pmatrix} 1
& r'O_K\\ 0 & 1\end{pmatrix}$.  Pick a generic element $g\in \Gamma $,
with    Bruhat   decomposition   of    the   form    $g=umwv$,   where
$m=\begin{pmatrix}  \alpha  & 0  \\0  &  \beta  \end{pmatrix}$ may  be
assumed to  be as  in the foregoing  paragraph. Then,  $\Gamma \supset
<^g(V^+),  V^+>$. Note  that  $u,v$ centralise  the  group $V^+$;  put
$V^-=^w(V^+)$.  One sees that $\Gamma \supset ^u<^m(V^-),V^+>=^u<
\begin{pmatrix}1 & 0\\ \alpha ^{-1}\beta r'O_K & 1\end{pmatrix}, 
\begin{pmatrix}1 & r'O_K\\ 0 & 1\end{pmatrix}>$. By the result on
SL(2)  over   CM  fields  (Proposition   \ref{CM'}),  $\Gamma  \supset
^u(\Delta)$  for   some  subgroup  $\Delta   $  of  finite   index  in
$SL_2(O_E)$, where $E=K(\alpha ^{-1}\beta)$  a CM extension of $K$. In
particular, there  exists an integer  $r''$ such that  $\Gamma \supset
^u(\theta ^{r''\Z})$. By  Proposition \ref{technical}, it follows that
$\Gamma $ is arithmetic.\\

\subsection{Groups of outer type A}

\subsubsection{\bf The Groups SU(h) over fields} In this subsection, $K$
is  a number  field,  $E/K$ a  quadratic  extension whose  non-trivial
Galois  automorphism   is  denoted  $\sigma$.    Let  $h:E^{n+1}\times
E^{n+1}\ra E$ denote a $\sigma$-hermitian form which is isotropic over
$K$, and write $h=\begin{pmatrix}0 & 1\\ 1 & 0\end{pmatrix}\oplus h_0$
where  $h_0$ is  {\bf anisotropic}  over  $K$.  Let  $G=SU(h)$ be  the
special  unitary  group  of   this  hermitian  form.   Then,  $K$-rank
($G$)=$1$. The positive roots are $\alpha$ and $2\alpha$.  Assume that
$\R$-rank    ($G_\infty$)     $\geq    2$.     Therefore,    $\R$-rank
$(SU(h_0))_\infty  \geq 1$.  The  arguments are  general when  $K$ has
infinitely many units or when $n$  is large (i.e. $n\geq 4$).  But for
small $n$ and small fields, the proofs become more complicated, and we
give ad hoc arguments. We thus have 5 cases to consider.\\

{\it Case 1.  $K$ has infinitely many units}.  Note that the $2\alpha$
root  space  is  one  dimensional.   Therefore, by  the  criterion  of
Proposition  \ref{2alpha},  every  arithmetic  group  is  virtually  3
generated. \\

{\it Case  2. $K$ is  $\Q$ or is  an imaginary quadratic  extension of
$\Q$ but $n\geq 4$}. Take $P$ (resp. $M$) to be the parabolic subgroup
of $G$  (resp. Levi  subgroup of $P$),  consisting of matrices  of the
form  $\begin{pmatrix}a &  *  & *  \\ 0  &  b &  *\\  0 &  0 &  \sigma
(a)^{-1}\end{pmatrix}$ (resp. $\begin{pmatrix}a & 0 & 0 \\ 0 & b & 0\\
0 & 0 &  \sigma (a)^{-1}\end{pmatrix}$).  Then, $(M\supset )M_0\supset
M_1=SU(h_0)$ since  the latter group is  semi-simple (because $n-2\geq
2$; we only  use the hypothesis that $n\geq  3$, so these observations
apply to the  next two cases as well) and  is non-compact at infinity,
and  therefore contains  an  arithmetic subgroup  as  a Zariski  dense
subgroup. Moreover, $SU(h_0)(\C)=SL_{n-1}(\C)$, and its representation
on the  lie algebra $\fu$ of  the unipotent radical of  $P$, is simply
$St\oplus   St^*\oplus  triv$.    Since   $M_1(\C)=SL_{n-1}(\C)$  with
$n-1\geq  3$, the  standard representation  is not  equivalent  to its
contragredient.   Thus  the diagonal  torus  $T_1$  of  $M_1$ has  one
dimensional  eigenspaces  in  $\fu$.   Hence arithmetic  subgroups  of
$G(K)$ are virtually three generated. \\
 
{\it Case 3.  $n=3$, either $K=\Q$ and $E/\Q$ is real quadratic or $K$
is    an   imaginary    quadratic   extension    of    $\Q$}.    Then,
$SU(h_0)(\C)=SL_2(\C)$, but the torus $T_1$  of the last case does not
have multiplicity  one in its  action on $\fu$. However,  observe that
$M_0$ of the last case contains in addition the torus $T_2$ consisting
of matrices $\begin{pmatrix}u & 0 & 0 & 0 \\ 0 & u^{-1} & 0 & 0 \\ 0 &
0 & u^{-1} & 0  \\ 0 & 0 & 0 & u\end{pmatrix}$ with  $u$ a unit in the
real  quadratic extension  $E$ ($E$  has infinitely  many  units). Put
$T_0=T_1T_2$.  Now, $T_0\subset M_0$ is a torus consisting of matrices
of the form $\begin{pmatrix}u & 0 & 0 &  0 \\ 0 & u^{-1}v & 0 & 0 \\ 0
& 0  & u^{-1}v^{-1} & 0  \\ 0 & 0  & 0 &  u\end{pmatrix}$ with $u,v\in
{\bf G}_m$, and  has ( as may be easily  seen) distinct eigenvalues in
$\fu  :  St\oplus St^*\oplus  triv=\fu$,  where  $\fu$  is as  in  the
previous case.  Thus, arithmetic  subgroups of $G=SU(h)$ are virtually
three generated.\\

{\it Case 4.  $K=\Q$, $n=3$  and $E/\Q$ is imaginary quadratic.} Then,
$U(h_0)$  is  not  contained  in  $M_0$  (of  course,  $SU(h_0)\subset
M_0$). We give an ad hoc argument in this particular case.\\

Write   $h_0=\begin{pmatrix}\lambda    _1   &   0\\    0   &   \lambda
_2\end{pmatrix}$ with  $\lambda _1,\lambda _2\in  \Q$ (every Hermitian
form in two variables is equivalent to one of this type). Now, $h=
\begin{pmatrix}0 & 1\\ 1 & 0\end{pmatrix}\oplus h_0$ is viewed a
{\bf hermitian} (with respect to $\sigma$) form from $E^4\times E^4\ra
E$. Consider $f=\begin{pmatrix}0 & 1\\ 1 & 0\end{pmatrix}\oplus
\begin{pmatrix}\lambda _1 & 0\\ 0 & \lambda _2\end{pmatrix}$ as a {\bf
quadratic } form on $\Q^4$.  Now, the $\Q$-group $SU(h)$ contains as a
$Q$-subgroup,  the group  $SO(f)$.  Since  $SU(h_0)_\infty=SU(1,1)$ is
non-compact (as we  have seen before, this follows  from the fact that
the  real   rank  of  $G_\infty$   is  $\geq  2$),  it   follows  that
$SO(f_0)_\infty=SO(1,1)$      is      also     non-compact.       Here
$f_0=\begin{pmatrix}\lambda _1  & 0\\ 0 &  \lambda _2\end{pmatrix}$ is
viewed as a quadratic  form. Consequently, the group $SO(f)(\R)\supset
SO(1,1)\times SO(1,1)$ and therefore has real rank $\geq 2$.\\

{\it Claim:  $SO(f)$ is a $\Q$-simple  group}. For, if  $SO(f)$ is not
$\Q$-simple,  (since it  is isogenous  to the  product $SL_2(\C)\times
SL_2(\C)$) then  it is isogenous  to $SL_2\times SL_2$  or $SL_2\times
SL_1(D)$  over $\Q$  (with $D$  a quaternionic  division  algebra over
$\Q$).    Now,   the   only   four  dimensional   representations   of
$SL_2(\C)\times  SL_2(\C)$ are  $St\otimes St$,  or $St\oplus  St$, or
$Triv\otimes Triv\oplus  Triv \otimes S^2(St)$,  or $Triv^2\oplus Triv
\otimes St$. Thus, if both the factors have to act non-trivially, then
the only  possible four dimensional representations  are $St\oplus St$
and $ST\otimes St$. But, $St\oplus  St$ does not have a quadratic form
invariant   under  $SL_2\times   SL_2$.   Thus,   the   only  possible
representation  (over  $\C$), of  $SL_2\times  SL_2$  onto $SO(4)$  is
$St\otimes St$.\\

It  follows that  the group  $SL_2\times SL_1(D)$  cannot have  a four
dimensional representation defined over $\Q$ with image $SO(f)$. Thus,
$SO(f)$ must be isogenous to  $SL_2\times SL_2$. But then, the isogeny
$SL_2\times  SL_2   \ra  SL(St\otimes  St)=SL(M_2(\Q))$   preserves  a
quadratic form (namely the determinant) over $\Q$, which has {\bf two}
$\Q$-hyperbolic planes  in it,  and therefore, cannot  be of  the form
$f=J\oplus f_0$ as above. This proves the claim. \\

Choose  an  element $\theta  \in  SO(h_0)$  of  infinite order.   Pick
non-trivial   elements  $u_0\in   (SO(f)\cap  U^+)(\Z)$   and  $v_0\in
U_{2\alpha}(\Z)$.   Then, the  group $<\theta,  u_0v_0>$  generated by
$\theta $ and $uv$ contains $\theta ^\Z$, and (since $\theta $ acts by
different  characters on  $Lie U_{2\alpha}$  and $U^+\cap  SO(f)$) the
unipotent  group  $V^+=V^+(r\Z)==SO(f)\cap U^+(r\Z)U_{2\alpha  }(r\Z)$
for some integer $r$. Let $\gamma  \in G(\Z)$ be an element in general
position  with  respect  to  $\theta   $  and  $u_0v_0$  as  in  Lemma
\ref{zariskidense}.  For  an integer  $r$, consider the  group $\Gamma
=<\theta  ^r, (u_0v_0)^r,  \gamma  ^r>$.  Then  $\Gamma  $ is  Zariski
dense. By  arguments similar to the  previous cases, it  to prove that
every arithmetic group in $G(\Z)$  is virtually three generated, it is
enough to  show that $\Gamma $  is arithmetic for every  $r$.  Pick an
element $g\in \Gamma $  with Bruhat decomposition $g=umwv$, say. Then,
there exists  an integer $r'$ such  that $u$ and  $v$ take $V^+(r'\Z)$
into  $V^+(r\Z)$ (since  commutator  of $u$  with  $V^+$ lands  inside
$U_{2\alpha}$).   Thus,   $\Gamma  \supset  ^u<^m(V^-),   V^+>$  where
$V^-=^w(V^+)$   as   before.     Hence,   we   get   $\Gamma   \supset
^u<U_{-2\alpha}(r\Z), V^+(r\Z)>$. \\

Consider the  group $<U_{-2\alpha}(r\Z), U_H^+(r\Z)>$.   An element in
$U_{-2\alpha }(r\Z)$  has the Bruhat  decomposition $u_1m_1wv_1$ where
$u_1,v_1\in U_{2\alpha }(\Q)$  {\it commute with} $U_H^+$.  Therefore,
$<U_{-2\alpha}(r\Z), U_H^+(r\Z)>$ contains
\[<^{u_1m_1wv_1}(U_H^+(r\Z)), U_H^+(r\Z)>=
^{u_1}<^{m_1}(U_H^-(r\Z)),  U_H^+(r\Z)>  \]  and the  latter  contains
$^{u_1}<U_H^-(r'\Z), U_H^+(r'\Z)>$  for some integer  $r'$.  Since the
latter group  is of  finite index in  $H(\Z)$ by \cite{V},  it follows
that     $<U_{-2\alpha}(r\Z),     U_H^+(r\Z)>\supset     ^{u_1}(\theta
^{r\Z})=\theta  ^{r\Z}$ for  some  integer $r$.   Therefore, from  the
foregoing   paragraph,   we   get   $\Gamma   \supset   ^{uu_1}(\theta
^{r\Z})=^u(\theta ^{r\Z})$.   Then, $\Gamma $  contains the commutator
$[^u(\theta  ^r),  \theta  ^r]$,  with  $u$  running  through  generic
elements of  $U^+$ whence, $\Gamma \supset U^+(r\Z)$  for some integer
$r$.  By \cite {V}, $\Gamma $ is arithmetic.\\

{\it  Case 5.  $n=2$,  $K$ is  either $\Q$  or an  imaginary quadratic
extension of $\Q$}.

We  can take  $h=\begin{pmatrix} 0  & 0  &  1\\0 &  1 &  0\\ 1  & 0  &
0\end{pmatrix}$ as  a Hermitian form over a  quadratic extension $E/K$
and $G=SU(h)$ over $K$.

If  $K=\Q$, and  $E$ is  imaginary quadratic,  then the  real  rank of
$SU(h)$ is one (the group of real points is $SU(2,1)$), and in Theorem
1 we have assumed that  $R-rank (G_\infty)\geq 2$.  Hence, $E$ is real
quadratic  and  has  infinitely  many  units.   If  $K$  is  imaginary
quadratic, then any quadratic extension $E$ of $K$ has infinitely many
units.   We  can  therefore  assume  that $E/K$  has  infinitely  many
units. \\

If  $P$ is  the parabolic  subgroup of  $G=SU(h)$ consisting  of upper
triangular matrices in $G$, then it follows from the conclusion of the
last paragraph, that $M_0(\C)=C^*$ since $M_0(O_K)$ contains the group
of matrices $h=\begin{pmatrix} u & 0 & 0\\0 & u^{-2} & 0\\ 0 & 0 & u
\end{pmatrix}$ where $u$ is a unit in $E$.  The action of
$M_0(\C)=C^*$ on the Lie algebra $\fu$ of the unipotent radical of $P$
is given by $\fu =\C(3)\oplus \C(-3)\oplus \C(0)$ where $\C(m)$ is the
one dimensional module over $M_0(\C)$  on which an element $z\in \C^*$
acts by  $z^m$.  Hence  $\fu$ is multiplicity  free for  the $M_0(\C)$
action, and we have proved Theorem 1 in this case. \\

\subsubsection{\bf The Groups SU(h) over division algebras} 
In  this  subsection,  $K$  is  a  number  field,  $E/K$  a  quadratic
extension, $D$ a central division  algebra over $E$ with an involution
$*$    of     the    second    kind,     degree    ($D$)=$d\geq    2$,
$k=[K:\Q]$. $h:D^{m+2}\times D^{m+2}\ra D$  is a $*$-hermitian form in
$m+2$ variables over $D$. $h$ is of the form
\[h=\begin{pmatrix}0 & 1\\ 1 & 0\end{pmatrix}\oplus h_m\] 
where $h_m$  is an  anisotropic hermitian form  in $m$  variables. The
special  unitary group  $G=SU(h)$  of  the hermitian  form  $h$ is  an
absolutely simple algebraic group over $K$, and under our assumptions,
$K$-rank($G$) $=1$.

{\it Case 1.  $D\otimes \R\neq {\bf H}\times\cdots \times {\bf H}$}.

Then,  the   group  $SL_1(D\otimes  \R)$   is  not  compact,   and  is
semi-simple. If $U^+=\{\begin{pmatrix}1 & * & *\\  0 & 1 & *\\ 0 & 0 &
1\end{pmatrix}\}$, $U_{2\alpha}=\{\begin{pmatrix}1  & 0 & x\\ 0  & 1 &
0\\ 0  & 0 & 1\end{pmatrix}:  x+x^*=0\}$ and $P$ is  the normaliser of
$U^+$ in $G$  (then $P$ is a parabolic subgroup of  $G$), there is the
obvious Levi subgroup  $M$ of $P$.  Since $SL_1(D)$  is non-compact at
infinity,    it     follows    that    $M_0\supset     M_1$,    where,
$M_1=R_{E/K}(SL_1(D))$.   Moreover, $M_1(\C)=SL_d(\C)\times SL_d(\C)$,
and as  a module over $M_1(\C)$,  the Lie algebra  $\fu _{2\alpha}$ of
$U_{2\alpha}$ is  $\C^d\otimes (\C^d)^*$.  Thus, the  weight spaces of
the  torus $T$= diagonal$\times  $ diagonal  of $SL_d\times  SL_d$, on
$\fu _{2\alpha}$  are all one dimensional. Hence,  there exist $m_0\in
M_1(O_K)$,  $u_0\in U_{2\alpha  }(O_K)$  such that  for every  integer
$r\geq 1$,  there exists  an integer $r_0$  with $<m_0^r,u_0^r>\supset
U_{2\alpha }(r_0O_K)$. \\

By  Lemma  \ref{zariskidense}, there  exists  an  element $\gamma  \in
G(O_K)$   such  that   for  any   integer  $r$,   the   group  $\gamma
=<m_0^r,u_0^r,\gamma ^r>$ is Zariski  dense in $G(K\otimes \C)$. As in
the  previous  sections,  it  suffices  to prove  that  $\Gamma  $  is
arithmetic.  Pick  $g=umwv\in \Gamma $.  Then, $\Gamma  $ contains for
some  integers   $r',r''$  the  groups  $^g(U_{2\alpha}(r'O_K))\supset
^u(U_{-2\alpha          }(r"O_K))$         as          well         as
$U_{2\alpha}(r'O_K)=^u(U_{2\alpha}(r'O_K)$. \\

Consider the group $H=SU(J,D)$,  where $J$ is the hyperbolic hermitian
(with  respect to  $*$)  form in  two  variables given  by the  matrix
$\begin{pmatrix}  0 &  1\\  1 &  0\end{pmatrix}$.   $H$ is  absolutely
simple over  $K$. Moreover, it  contains as a $K$-subgroup,  the group
$M_H=R_{E/K}(GL_1(D))$ the embedding given by $g\mapsto
\begin{pmatrix} g & 0\\ 0 & (g^*)^{-1}\end{pmatrix}$.  Now,
$M_H(K\otimes  \R)=GL_1(D\otimes  \R)\supset  R^*\times  SL_1(D\otimes
\R)$.   Since $SL_1(D\otimes  \R)$ is  not compact  by  assumption, it
follows that  $M_H(K\otimes \R)$ has  real rank $\geq 2$.   The groups
$U_{\pm 2\alpha}$ are maximal opposing unipotent subgroups of $H$, and
hence            by           \cite{V},            the           group
$<U_{2\alpha}(r'O_K),U_{-2\alpha}(r''O_K)>$ is  an arithmetic subgroup
of $H(O_K)$. Therefore,  we get from the last  paragraph, that $\Gamma
\supset  ^u(\Delta')$ for  some subgroup  $\Delta' \subset  H(O_K)$ of
finite index, which implies  that $\Gamma \supset ^u(\Delta)$ for some
subgroup $\Delta $ of finite index in $SL_1(O_D)$ for some order $O_D$
in  $D$.   Since  $SL_1(O_D)$  contains  elements which  do  not  have
eigenvalue 1 in their action on $Lie U^+$, it follows from Proposition
\ref{technical} that $\Gamma $ is arithmetic. \\

{\it  Case  2. $D\otimes  \R={\bf  H}\times\cdots\times  {\bf H}$  and
$m\geq 2$}.
 
Then  $E$  is  totally real  (and  so  is  $K$),  and  $D$ must  be  a
quaternionic division  algebra over $E$.   Moreover, $SU(h_m)(K\otimes
\R)=\{g\in SL_m(D\otimes \R): g^*h_mg=h_m\}=\{g=(g_1,g_2)\in SL_m({\bf
H})^k               \times               SL_m({\bf              H})^k:
(g_2^{\iota},g_1^{\iota}(h_m,h_m)(g_1,g_2)=(h_m,h_m)\}$  where $\iota$
is  the  standard  involution  on  ${\bf  H}$  induced  to  $SL_m({\bf
H})^k$. Thus, $SL_m(K\otimes \R)$  is isomorphic to $SL_m({\bf H})^k$.
Since  $m\geq  2$, the  group  $SL_m({\bf  H})^k$  is semi-simple  and
non-compact, and  contains a  Zariski dense set  of integral  points ,
which  are  $SU(h_m)(O_K)=M_1(O_K)$.   Take  $P$ to  be  the  standard
parabolic    subgroup    of    $G=SU(h)$.     Hence    $M_0(\C)\supset
M_1(\C)=SL_{2m}(\C)^k$.  As  a module over $M_1(\C)$,  the Lie algebra
$Lie      U^+(\C)=[\C^2\otimes     (C^{2m})^*\oplus     \C^{2m}\otimes
(\C^2)^*\oplus triv  ^4]^k$.  Choose  a generic toral  element $m_0\in
M_1(O_K)$,  and  an  element  $u_0=u_1u_1\in  U^+(O_K)$  with  $u_1\in
Exp(\fg_\alpha)$ and $u_2\in  Exp(\fg _{2\alpha})$.  Choose an element
$\gamma  \in  G(O_K)$ of  infinite  order,  in  general position  with
respect to  $u$ and $m$  (Lemma \ref{zariskidense}).  Then,  for every
integer $r$,  the group $\Gamma =<u_0^r,\gamma ^r,  m_0^r>$ is Zariski
dense. \\

Let $\Delta\subset U^+$ be the group generated by $^{m_0^{jr}}(u_0^r):
j\in \Z$ and $Log :U^+\ra \fu$  the log mapping. Then, $log (\Delta )$
contains  elements of  the form  $v_1,\cdots,v_N$ with  each  $v_i$ an
eigenvector for  $m_0\in \prod SL_{2m}(\C)=M_1(\C)$.   For the generic
toral  element $m_0$,  the  number of  distinct  {\it eigenvalues}  on
$V_1=(\C^2)^*\otimes   \C^{2m}\oplus  \cdots   \oplus  (\C^2)^*\otimes
\C^{2m}$ (the direct sum taken  $k$ times) is $2mk$. Fix corresponding
eigenvectors  $v_1^i,\cdots, v_{2m}^i~(1\leq i\leq  k$ in  $V_1$. Pick
similarly, $2mk$ eigenvectors  $(v_1^i)^*,\cdots (v_{2m}^I)^* ~ (1\leq
i\leq  k$  in $V_1^*=  (\C^2)\otimes  (\C^{2m})^*\oplus \cdots  \oplus
(\C^2)\otimes  (\C^{2m})^*$  (the  direct  sum taken  $k$  times)  for
$m_0$. The trivial  $M_1(\C)$ module $\fg _{2\alpha}$ is  the $k$ fold
direct sum  of $M_2(\C)$  with itself. Denote  the $i$th  component of
this  direct sum $M_2(\C)_i$  ($1\leq i\leq  k$). By  general position
arguments (since the toral element  $m_0$ is generic) it can be proved
that  for each  $i$,  the  $2m-1 (\geq  3)$  vectors $  v_1^i(v_2^i)^*
-v_2^i(v_1^i)^*,    \cdots,    v_1^i(v_{2m}^i)^*   -v_{2m}^i(v_1^i)^*$
together with the vector $v_2^i(v_3^i)^* -v_3^i(v_2^i)^*$, span all of
the  $i$-th component  $M_2(\C)_i$.   We choose  $u_1$  such that  the
element  $log  (u_1)$  has  non-zero  projections  into  each  of  the
eigenspaces of $m_0$ in  $V_1\oplus V_1^*$, and its projections $v_\mu
^i,(v_\mu ^i)^*$ are as in the foregoing. \\

Therefore,  $\Gamma \supset  Exp (\Delta  )\supset U_{2\alpha}(r'O_K)$
for some  integer $r'$. To prove  Theorem 1 in this  case, by standard
arguments, it is  enough to prove that $\Gamma$  is arithmetic. Take a
generic  element $g=umwv\in \Gamma  $.  Then,  $\Gamma $  contains the
group   $<^g   (U_{-2\alpha}(r'O_K)),   U_{2\alpha}(r'O_K)u_1^{r'\Z}>$
(recall that $u_1\in Exp(\fg _\alpha)$).  Thus, for some other integer
(denoted  again  by  $r'$   to  save  notation),  $\Gamma  $  contains
$^u<U_{-2\alpha}(r'O_K), u_1^{r'\Z}U_{2\alpha}(r'O_K)>$. \\

Let us  view $h_0=\begin{pmatrix}  0 &  0 &  1\\0 & 1  & 0\\  1 &  0 &
0\end{pmatrix}$ as  a Hermitian form for  $E/K$.  Set $H=SU(h_0)\simeq
SU(2,1)$. This is  an algebraic group over $K$,  and has corresponding
upper triangular  unipotent group $U_H^+$. The Lie  algebra spanned by
$Elog(u_1)$ and  $Elog (^w(u_1)$  is easily seen  to be  isomorphic to
that  of $H$  with $Lie(U_H^+)=Elog  (u_1)\oplus [Elog  u_1,Elog u_1]$
(the square  bracket denotes the  commutator). From the  conclusion of
the  last paragraph, we  get $\Gamma  \supset <^u(U_{-2\alpha}(r'O_K),
u_1U_{2\alpha}(r'O_K)>$.   By  \cite {V},  the  latter group  contains
$^u(SU(2,1)(r'O_K))$. Hence $\Gamma $ contains $^u(SU(2,1)(r'O_K))$ as
$g=umwv$   varies,   and  for   some   fixed  $g'=u'm'wv'$,   contains
$^{u'}(SU(2,1)(r'O_K))$ as well.\\

The toral element $h\in  SU(2,1)$ of the form $h=\begin{pmatrix}\theta
& 0 & 0\\ 0 & \theta ^{-2} & 0\\ 0 & 0 & \theta \end{pmatrix}$ acts on
the root space $\fg _\alpha$ by the eigenvalues $theta, \cdots, \theta
$ and  $\theta ^3$  (as may  be easily seen).   Therefore, $h$  has no
fixed vectors in $\fg _\alpha$.  Now, by the last paragraph, $\Gamma $
contains  the group  $^u (h)$  ($u$ generic).   Hence,  by Proposition
\ref{technical}, $\Gamma $ is arithmetic. \\

{\it Case 3. $D\otimes \R={\bf H}^{2k}$ and $m\leq 1$, but $k\geq 2$}.

Again, $E$ and $K$ are  totally real. $G=SU(h)$ with $G(K\otimes \R)=$
$SL_3({\bf H})^k$ if $m=1$ and $SL_2({\bf H})^k$ if $m=0$. Fix $\theta
\in O_K^*$ such that $\Z[\theta ^r]$  is a subgroup of finite index in
$O_K$  for all  $r\neq 0$  (Lemma \ref{exist}).   Let $\alpha  $  be a
totally  positive element such  that $E=K(\sqrt{\alpha})$.   Denote by
$t(\theta )$ (resp.  $u_{+}$)  the matrix $\begin{pmatrix}\theta & 0 &
0\\0   &   1   &   0\\    0   &   0   &   \theta   ^{-1}\end{pmatrix}$
(resp. $\begin{pmatrix}1 & 0 & \sqrt{\alpha}\\0 & 1 & 0\\ 0 & 0 & 1
\end{pmatrix}$) if $m=1$ and the
matrix  $\begin{pmatrix}\theta &  0\\ 0  &  \theta ^{-1}\end{pmatrix}$
(resp.   $\begin{pmatrix}1 & \sqrt{\alpha}\\  0 &  1\end{pmatrix}$) if
$m=0$. By  the choice of $\theta  $, the group  $<\theta ^r, u_{+}^r>$
contains, for  every $r$, $u_{+}^{r'O_K}$ for some  integer $r'$. Pick
an element $\gamma  \in G(O_K)$ in general position  as in Proposition
\ref{zariskidense}.     Then   for    every    $r\neq   0$,    $\Gamma
=<t^r,u_{+}^r,\gamma   ^r>\subset   G(O_K)$   is  Zariski   dense   in
$G(K\otimes \C)$. Pick a generic element $g=umwv\in \Gamma $. Then,
\[\Gamma \supset 
<^g(u_{+}^{rO_K}),  u_{+}^{rO_K}>=^u<^m((u_-)^{rO_K}, u_{+}^{rO_K}>.\]
The element $m$ is of the form  $\begin{pmatrix}a & 0 & 0\\0 & b & 0\\
0 &  0 &  (a^*) ^{-1}\end{pmatrix}$ with  $b\in SU(h_m)$ if  $m=1$ and
$\begin{pmatrix}a & 0\\ 0 & (a^*)^{-1}\end{pmatrix}$ if $m=0$ for some
$a\in D^*$. Hence $^m(u_-^{rO_K})$=$\begin{pmatrix} 1 & 0 & 0\\0 & 1 &
0\\    \sqrt{\alpha}(aa^*)^{-1}r'O_K    &    0   &    1\end{pmatrix}$,
$u_{+}^{rO_K}=\begin{pmatrix}1 & 0 & r'O_K\\0 & 1 & 0\\ 0 & 0 & 1
\end{pmatrix}$ if $m=1$ 
(and   $^m(u_-^{r'O_K})=\begin{pmatrix}1  &  0\\   (aa^*)^{-1}r'O_K  &
1\end{pmatrix}$,   $u_{+}^{r'O_K}=\begin{pmatrix}1  &   r'O_K\\   0  &
1\end{pmatrix}$  if $m=0$).   These two  groups $^m(u_-)$  and $u_{+}$
generate $SL_2$ over $K(c)$ if $m=1$ and $SL_2$ over $K$if $m=0$. \\

The  element  $c=aa^*\in  D^*$  has  its reduced  norm  and  trace  in
$E$. But, in  fact, $Tr (c)$ and  $N(c)$ lie in $K$ itself,  as may be
easily seen. Now,  $c$ being in the quaternionic  division algebra $D$
over  $E$  with  $D\otimes   \R={\bf  H}^{2k}$,  generates  a  totally
imaginary  quadratic  extension  over  the totally  real  $E$.   Hence
$K(c)/K$ is also totally  imaginary quadratic extension.  By the SL(2)
result Proposition \ref{CM'} ($K$ is  a totally real number field with
infinitely many units), we get $<^m(u_-^{r'O_K}), u_{+}^{r'O_K}>$ is a
subgroup of finite index  in $SL_2(O_{K(c)})$ if $m=1$ and $SL_2(O_K)$
if $m=0$.  In particular, the group $<^m(u_-^{r'O_K}), u_{+}^{r'O_K}>$
contains  the group  $t^{r'\Z}=t(\theta )^{r''\Z}$  for  some $r''\neq
0$.\\
 
Thus, $\Gamma $ contains the group $^u(t^{r''\Z})$, $u$ is generic and
$t$ does  not have  eigenvalue one  in its action  on the  Lie algebra
$Lie~U^+$.   Therefore, by Proposition  \ref{technical}, $\Gamma  $ is
arithmetic. \\

{\it  Case 4.   $D\otimes  \R={\bf H}^{2k}$,  $m=1$  and $k=1$}  (i.e.
$K=\Q$).  In  this case, we will explicitly  exhibit elements $u_{+}$,
$u_-$ and  $t$ in $G(O_K)$  such that for  every $r\neq 0$,  the group
$\Gamma =<u_+^r,u_-^r,t^r>$ is arithmetic.   This will prove Theorem 1
in  this case. Since  $D\otimes \R$  is a  product of  the Hamiltonian
quaternions ${\bf H}$, it follows that $E/\Q$ is real quadratic. Fix a
generic element $a\in  D^*$. Then, the element $aa^*$  generates as in
the last  case, an imaginary  quadratic extension over $\Q$.   Pick an
element  $t_2\in \Q(aa^*)\setminus  \Q$ such  that $t_2^2\in  \Q$. Now
choose $t_1\in  D$ such  that $t_1^2\in E$  but does not  commute with
$t_2$.  Write $E=\Q(\sqrt{z})$  where $z \in \Q$ is  positive.  Pick a
unit $\theta \in O_E^*$ of infinite order. \\

Write
\[u_+=\begin{pmatrix}1 & 1 & -\frac{1}{2}\\0 & 1 & -1\\ 0 & 0 & 1
\end{pmatrix}\begin{pmatrix}1 & 0 & \sqrt{\alpha}
\\0  & 1 &  0\\ 0  & 0  & 1\end{pmatrix},  u_-=\begin{pmatrix}1 &  0 &
0\\t_1 & 1 & 0\\ -\frac{t_1^2}{2} & -t_1 & 1
\end{pmatrix}\begin{pmatrix}1 & 0 & 0\\0 & 1 & 0\\
t_2\sqrt{\alpha} & 0 & 1\end{pmatrix}, \] and $t=\begin{pmatrix}\theta
& 0 & 0\\0 & \theta ^{-2} & 0\\ 0 & 0 & \theta\end{pmatrix}$. Now, the
group  $H=SU(2,1)$ for  the extension  $E/\Q$ embeds  in $G$  with the
corresponding group  of upper and lower  triangular unipotent matrices
$U_H^{\pm}$.  By Proposition \ref{SU(2,1)} applied to this SU(2,1), we
see  that  $\Gamma  $  contains  $U_H^{\pm}(r'O_{\Q(t_2)})$  for  some
integer $r'$. \\

In  particular,  $\Gamma  $  contains,  for some  $r'$,  the  subgroup
$U_H(r'O_{E\otimes  F})$   where  $E\otimes   F$  is  a   {\it  field}
($F=\Q(t_2)$ is  imaginary quadratic,  and $E/\Q$ is  real quadratic).
Taking  commutators  with  $u_-$,  we obtain  that$\Gamma  $  contains
elements  $v_-\in U_{-2\alpha}$ of  the form  $\begin{pmatrix}1 &  0 &
0\\0 &  1 & 0\\ r'x  & 0 &  1\end{pmatrix}$ with $x$ in  the subgroups
$\Z$, $t_1\Z$,  $t_2\Z$ and $t_1t_2\Z$  of $O_D$. However, the  sum of
these subgroups is a subgroup of finite index in $O_D$. Therefore, for
some  other $r'\neq  0$,  we get  $\Gamma \supset  U_{-2\alpha}(r'\Z)$
(i.e.  $\Gamma $ intersects $U_{-2\alpha}$ in an arithmetic subgroup).
Then, by Proposition \ref{highestroot}, $\Gamma $ is arithmetic. \\

{\it  Case  5.   $D\otimes  \R=  {\bf H}^{2k}$,  $m=0$,  $k=1$}  (i.e.
$K=\Q$).   Then, $G(\R)=SU(h)(\R)=SL_2({\bf H})$.   Therefore, $G(\R)$
has real rank one, and in Theorem 1, this is excluded. \\

\subsection{\bf Groups of type B} 

$G=SO(f)$  with $f$ a  non-degenerate quadratic  form in  $2l+1\geq 5$
variables  over a  number  field $K$.   $f$  is the  direct  sum of  a
hyperbolic form and an anisotropic form in $2l-1$ variables:
\[f=\begin{pmatrix} 0 & 1\\ 1 & 0\end{pmatrix}\oplus f_m \] 
with  $m=2l-1\geq 3$.   Assume that  $\R$-rank ($G_\infty$)  $\geq 2$,
where $G_\infty=G(K\otimes \R)$. Take $P$ to be the parabolic subgroup
$P=\{\begin{pmatrix}a   &  *   &   *\\0  &   b   &  *\\   0   &  0   &
a^{-1}\end{pmatrix}\in G: a\in GL_1/K, b\in SO(f_m)\}$.  The unipotent
radical   $U^+$   of   $P$   consists   of  matrices   of   the   form
$\begin{pmatrix}1 & x & -\frac{\sum x_i^2}{2}\\0 & 1_m & -^tx \\ 0 & 0
&  1\end{pmatrix}$ with  $1_m$ the  $m\times m$  identity  matrix, and
$x\in {\bf  A}^{2l-1}$, the affine  $2l-1$-space over $K$.   Denote by
$\fu $ the  Lie algebra of $U^+$. Let $U^-$ be  the transpose of $U^+$
(it  lies in  $G$).  Let  $M$ be  the Levi  subgroup of  $P$  given by
$M=\{\begin{pmatrix}a   &  0   &   0\\0  &   b   &  0\\   0   &  0   &
a^{-1}\end{pmatrix}\in   P:  a\in   GL_1/K,   b\in  SO(f_m)\}$.    Put
$H=SO(f_m)$.\\

{\it  Case  1.   $H_\infty=H(K\otimes  \R)$  is  non-compact}.   Then,
$H_\infty$  is  a non-compact  semi-simple  group,  hence $H(O_K)$  is
Zariski dense in $H(K\otimes \C)$.   Therefore, $M_0 \supset H$.  As a
module  over  $H(\C)=SO(2l-1,\C)$,   $\fu  (\C)=St=\C^{2l-1}$  is  the
standard representation,  and the maximal  torus $T_H$ of  $H(\C)$ has
distinct eigenvalues. Hence by Proposition \ref{multone}, Theorem 1 is
true for $G=SO(f)$ in this case. \\

{\it  Case 2.  $H_\infty=H(K\otimes  \R)$ is  compact}.  Then,  $K$ is
totally   real,   and   $(2\leq  )\R$-rank   ($G_\infty$)   =$\R$-rank
($GL_1(K\times  \R)$)=$[K:\Q]$, therefore $K\neq  \Q$.  Now,  $M_0$ is
rather  small. $M(O_K)$  is commensurate  to  $GL_1(O_K)=O_K^*$, hence
$M_0=  \{\begin{pmatrix}a   &  0  &   0\\0  &  1   &  0\\  0  &   0  &
a^{-1}\end{pmatrix}\in G:  a\in GL_1/K\}$.  In this case,  we will use
the  fact  that $SO(f)$  contains  many  $PSL_2(E)$  for totally  {\it
imaginary}  quadratic  extensions of  the  totally  real number  field
$K$. To see this, we first  prove a lemma.  Write the anisotropic form
$f_m$ as a direct sum $f_m=\phi\oplus \phi '$ with $\phi $ a quadratic
form  in {\bf  two} variables;  write $\phi  =1 \oplus  \lambda)$ with
$\lambda \in K$. Form the  quadratic forms $Q=\begin{pmatrix}0 & 1\\ 1
& 0\end{pmatrix} \oplus \phi$, and $Q_1=
\begin{pmatrix}0 & 1\\ 1 & 0\end{pmatrix} \oplus 1$. 
Then,    for    any    archimedean    completion   $K_v$    of    $K$,
$SO(Q)(K_v)=SO(3,1)(\R)\simeq  SL_2(\C)$.   Let  $Spin(Q)$ denote  the
simply connected two sheeted cover of $SO(f)$.

\begin{lemma} There exists a totally imaginary quadratic extension
$E/K$   such   that  $Spin(Q)$   is   $K$-isomorphic   to  the   group
$R_{E/K}(SL_2)$  where  $R_{E/K}$  denotes  the  Weil  restriction  of
scalars.
\end{lemma}

\begin{proof} Clearly, $SO(Q)$ is $K$-simple. Hence
$Spin  (Q)=R_{E/K}(H_0)$  with   $H_0$  an  absolutely  simple  simply
connected group over $E$, for some extension $E/K$, say of degree $d$.
Since  $SO(Q)$ is  isotropic over  $K$, so  is $H_0$  over  $E$. Since
$dim(SO(Q)/K)=6$,  one  sees that  $dim  (H_0)=\frac{6}{d}$. But  $dim
(H_0)\geq 3$  since it is  absolutely simple, hence $d\leq  2$.  Since
$Q$ is a form in four  variables, $Spin (Q)$ is not absolutely simple.
Therefore, $d=2$ (i.e.  $E/K$ is a quadratic extension), and $H_0$ has
dimension  $3$ (and  is isotropic  over $E$).   Therefore, $H_0=SL_2$.
Since     $SO(Q)(K_v)=PSL_2(\C)$,     it     follows    that     $Spin
(K_v)=SL_2(\C)=H_0(E\otimes K_v)$, for  every archimedean (hence real)
completion of $K$. Hence $E$ is a totally imaginary.
\end{proof}

The  inclusion of the  quadratic spaces  $Q_1$ and  $Q$ in  $f$ induce
inclusions  of  $SO(Q_1)$  and   $SO(Q)$  into  $SO(f)$  defined  over
$K$. They  further induce corresponding inclusions  (defined over $K$)
of  the  group of  unipotent  upper  (and  lower) triangular  matrices
$U_{Q_1}^{\pm}$ and  $U_Q^{\pm}$ into  the group $U^{\pm}$  defined at
the  beginning  of  this  subsection.  Let  $v\in  U^+_Q(O_K)\setminus
U_{Q_1}(O_K)$ ($U^-_{Q_1}$ is one dimensional). Then the $SL_2$ result
(Proposition \ref{CM'}) shows that $<v^{rO_K}, U^-(rO_K)>$ generates a
subgroup  of finite index  in $SO(Q)(O_K)$  (which is  commensurate to
$SL_2(O_E)$). \\

Let  $H=SO(Q_1)\subset SO(f)$,  $U^+_H=U^+\cap H$  be as  in  the last
paragraph. Let $t=\begin{pmatrix}\theta & 0 & 0\\0 & 1_m & 0\\ 0 & 0 &
\theta  ^{-1}\end{pmatrix}$ be  in $M(O_K)$,  $\theta \in  O_K^*$ such
that $\Z[\theta ^r]$ of finite index in $O_K$ (Lemma \ref{exist}). Fix
$u_+\in U^+_H(O_K)$, $u_+\neq 1$. Let $\gamma $ be in general position
with respect to $t $ and $u_+$. Write, for an integer $r\neq 0$,
\[\Gamma =<u_+^r, t^r,\gamma ^r>.\]
Then $\Gamma $ is Zariski  dense in $G(K\otimes \C)$. Moreover, by the
assumptions   on   $\theta$,   $\Gamma   $   contains   the   subgroup
$V^+(r')=u_+^{r'O_K}$for   some  $r'$.    Define   $V^-(r')$  as   the
$w$-conjugate of $V^+(r')$. Pick  a generic element $g=um'wv\in \Gamma
$.       Then,       $\Gamma      $      contains       the      group
$<^g(V^+(r')),V^+>=^{um'w}(V^+(r''),  V^+(r')>$ for  some  $r''$.  The
latter group contains $^u<^{m'}(V^-(r'')), V^(r'')>$ (replace $r''$ by
a larger $r''$ if necessary). \\

If $log  u_+=X\in Lie U^+\simeq K^m$,  then for the  generic $m'$, the
vectors  $X$ and  $^{m'}(X)$ span  a two  dimensional subspace  of the
anisotropic  quadratic space  $(K^m,f_m)$.  Write  the  restriction of
$f_m$ to $W$ as $\mu \phi$ for some $\mu \in K$, and $\phi $ as in the
Lemma  above, with  $\phi  (X,X)=1$, say.   Then,  $\Gamma $  contains
$^u<^{m'}(exp  (r''O_KX),   exp  (r''O_KX)>$,  which   by  Proposition
\ref{CM'},  contains $^u(\Delta  )$  for some  subgroup  $\Delta $  of
finite  index  in $SO(Q)(O_K)$,  where  $Q$  is  the four  dimensional
quadratic form as  in the Lemma.  Now, $\Delta  $ contains $t^{r_0\Z}$
for  some  $r_0$.  Hence  $\Gamma  $  contains  $^u(t^r)$, with  $t\in
M_0(O_K)\cap \Gamma  $.  By Proposition \ref{technical},  $\Gamma $ is
arithmetic.  This  proves Theorem  1 for $K$-rank  one groups  of type
B. \\

\subsection{Groups of type C} The groups of type C are $Sp_{2n}$ over
$K$  (which does  not have  $K$-rank 1),  and certain  special unitary
groups over  quaternionic division algebras.  In the  case of $K$-rank
one groups, we need only consider the groups of the latter kind. Thus,
let $D$ be a quaternionic central division algebra over $K$, $\sigma :
D\ra D$ an involution of the  {\it first} kind, such that the space of
$\sigma $ invariants in  $D$ is precisely $K$: $D^{\sigma}=K$. Suppose
$h:D^n\times D^n\ra D$ is a  $\sigma$-hermitian form which is a sum of
a hyperbolic form in two variables and an anisotropic form:
\[h=
\begin{pmatrix} 0 & 1\\ 1 & 0\end{pmatrix}
\oplus  h_{n-2}\]  with $h_{n-2}$  an  anisotropic  hermitian form  on
$D^{n-2}$.  The subgroup $P$ of $G$ consisting of matrices of the form
$
\begin{pmatrix}g & 0 &
  0\\0 & h & 0\\0 & 0 & (g^{\sigma})^{-1}\end{pmatrix}
\begin{pmatrix}1 & z &
  w\\0  & 1_{n-2}  &  0\\0 &  -^tz  & 1\end{pmatrix}$  is a  parabolic
subgroup   with  unipotent  radical   $U^+$  consisting   of  matrices
$\begin{pmatrix}1 & z & w\\0 & 1_{n-2} & 0\\0 & -^tz & 1\end{pmatrix}$
with $w+w^{\sigma}=0$. The commutator of $U^+$ is $U_{2\alpha}$ is the
set of matrices $\begin{pmatrix}1  & 0 & w\\ 0 & 1_{n-2}  & 0\\0 & 0 &
1\end{pmatrix}$  with  $w+w^{\sigma}=0$   having  dimension  $3$  over
$K$. Then $M$  is the Levi subgroup of $P$, with  elements of the form
$\begin{pmatrix}g   &    0   &   0\\0   &    h   &   0\\0    &   0   &
(g^{\sigma})^{-1}\end{pmatrix}$.\\

{\it Case 1.  $D\otimes \R\neq {\bf H}\times \cdots \times {\bf H}$}.

Then,  $SL_1(D\otimes   \R)$  is  a   non-compact  semi-simple  group.
Therefore, $M_0$  contains $SL_1(D)$, embedded as the  subgroup of $M$
of matrices of the form $\begin{pmatrix}g  & 0 & 0\\0 & 1_{n-2} & 0\\0
& 0  & (g^{\sigma})^{-1}\end{pmatrix}$ with $g\in  SL_1(D)$. Note that
for   any  embedding   of  $K$   in  $\C$,   we   have  $SL_1(D\otimes
_K\C)=SL_2(\C)$. As  a representation  of $SL_2(\C)$, the  module $Lie
U_{2\alpha}$  is $Sym  ^2(\C^2)$ (since  the space  $w=-w^{\sigma}$ is
3-dimensional) which  is multiplicity free  for the diagonal  torus in
$SL_2(\C)$. Therefore, there exists an $m_0\in SL_1(O_D)$, and $u_0\in
U_{2\alpha}((O_K)$  such  that the  group  generated  by the  elements
$^{m_0^j}(u_0): j\in \Z$ has finite index in $U_{2\alpha}(O_K)$. \\

Choose an element $\gamma \in G(O_K)$ in general position with respect
to $m_0,u_0$  as in Lemma  \ref{zariskidense}.  Write, for  an integer
$r\neq  0$, $\Gamma  =<m_0^r,u_0^r, \gamma  ^r>$. Then,  $\Gamma  $ is
Zariski dense in  $G(K\otimes \C)$.  By Proposition \ref{highestroot},
$\Gamma $ is arithmetic. This proves Theorem 1 in this case. \\

{\it Case 2.  $D\otimes \R={\bf H}^k$, $k=[K:\Q]\geq 2$}.  Then $K$ is
totally real, and contains an element $\theta \in O_K^*$ such that the
ring generated by $\theta $ has  finite index in the integers $O_K$ of
$K$.   Pick a  non-trivial element  $u_+\in U_{2\alpha}(O_K)$  and let
$u_-$  denotes its  conjugate  by  the Weyl  group  element $w$.   Let
$t=\begin{pmatrix}\theta &  0 &  0\\0 &  1_{n-2} & 0\\0  & 0  & \theta
^{-1}\end{pmatrix}\in   G(O_K)$.     For   $r\neq   0$,    the   group
$<t^r,u_+^r>\supset   u_+^{r'O_K}=V^+$   for   some  integer   $R'\neq
0$. Choose  $\gamma \in  G(O_K)$ in general  position with  respect to
$t,u_0$.   Then, $\Gamma  =<u_+^r,t^r,  \gamma ^r>$  is Zariski  dense
(Lemma \ref{zariskidense}). Then, $V^+\subset \Gamma $. Pick a generic
element $g=umwv\in  \Gamma $.  Then,  $\Gamma $ contains  the subgroup
$<^g(V^+),V^+>\supset  ^u<^m(u_+^{r''O_K}),   u_+^{r'O_K}>$  for  some
other integer $r''$. If $u_+=\begin{pmatrix}1 &  0 & w\\0 & 1 & 0\\0 &
0 & 1\end{pmatrix} $ then $^m(u_+)=\begin{pmatrix}  1 & 0 & 0\\0 & 1 &
0\\(a^{\sigma})^{-1}wa^{-1}    &    0    &    1\end{pmatrix}$    where
$m=\begin{pmatrix}a   &   0   &   0\\0    &   1   &   0\\0   &   0   &
(a^{\sigma})^{-1}\end{pmatrix}$.   Since $m$  is generic,  the element
$\xi  =(a^{\sigma  })^{-1}wa^{-1}w^{-1}\in  D$ generates  a  quadratic
(totally imaginary, by  the assumption on $D$ in  this case) extension
of   $K$.    Therefore,    by   Proposition   \ref{CM'},   the   group
$<^m(u_-^{r'O_K}),u_+^{r'O_K}>$   is   an   arithmetic   subgroup   of
$SL_2(K(\xi))$. In particular,  $\Gamma $ contains $^u(t^{r_0\Z})$ for
some  integer $r_0$.   The action  of $t$  on $Lie  U^+$ has  no fixed
vectors. By Proposition \ref{technical}, $\Gamma $ is arithmetic. \\

{\it Case 3. $D\otimes \R={\bf H}$, $k=1$ (i.e.  $K=\Q$)}.

Let $H=SU(h_{n-2})$.  Since $\R$-rank  ($G$) $\geq 2$, it follows that
$\R$-rank   ($H$)   $\geq  1$.    Thus,   $n-2\geq   2$.   But   then,
$H(\R)=SU(h_{n-2},{\bf  H})$ is  isotropic  if and  only if  $h_{n-2}$
represents a zero over $\R$.   Since $n-2\geq 2$, and a hermitian form
in $\geq 2$ variables over  a quaternionic algebra (with respect to an
involution of the first kind  whose fixed points are of dimension one)
represents a zero  over $\Q_p$ for every prime $p$,  it follows by the
Hasse principle (  see Ch 6, section (6.6),  Claim (6.2) of \cite{PR})
that $h_{n-2}$ represents a a zero over $\Q$ as well, whence $\Q$-rank
($H$) $\geq  1$ and $\Q$-rank ($G$)  $\geq 2$; this case  is not under
consideration in this section. \\

\subsection{Classical groups of type D} 

{\it   Case   1.   $G=SO(f)$}.    Here,   $f=J\oplus  f_{2n-2}$   with
$J=\begin{pmatrix} 0  & 1\\ 1  & 0\end{pmatrix}$ being  the hyperbolic
form  on $K^^2$, $f_{2n-2}$  an anisotropic  quadratic form  in $2n-2$
variables over $K$,  and $n\geq 4$ (i.e. $n-1\geq  3$).  Now, the real
rank of $SO(f)(K\otimes \R)$ is  $\geq 2$.  The argument for groups of
type B applies without change.  \\

We now assume that $G=SU_n(h,D)$.  Here, $D$ is a quaternionic central
division algebra  over $K$  with an involution  $\sigma$ of  the first
kind such that  the dimension of the set  of fixed points $D^{\sigma}$
is three (in the symplectic  case, this dimension was one). $h=J\oplus
h_{n-2}$, where $ J=\begin{pmatrix} 0 & 1\\ 1 & 0\end{pmatrix}$ is the
hyperbolic  form on $D^2$  and $h_{n-2}$  is an  anisotropic hermitian
form. Let  $P$, $U^+$ and  $M$ be as  in the symplectic case  The only
change  from that  case  is  that the  involution  $\sigma$ has  three
dimensional  fixed  space,   hence  $U_{2\alpha}$  which  consists  of
matrices of the form $\begin{pmatrix}1 & 0 & w\\0 & 1_{n-2} & 0\\0 & 0
& 1\end{pmatrix}$ with $w+w^{\sigma}$ is one dimensional over $K$.

{\it Case  2.  $K$  has infinitely many  units}. Then,  by Proposition
\ref{2alpha}, Theorem 1 holds in this case. \\

{\it Case 3.  $K$ is  an imaginary quadratic extension of $\Q$}. Then,
$SL_1(D\otimes    \R)=SL_2(\C)$.     Moreover,   $SU(h_{n-2})(K\otimes
\R)=SO_{2n-4}(\C)$. Note that $n\geq 4$. Therefore, $SO_{2n-4}(\C)$ is
a  semi-simple   group.   Hence,  $M_0(K\otimes   \R)=  M_0(\C)\supset
SL_2(\C)\times SO_{2n-4}(\C)$.  Then, the  product of the diagonals in
the latter group has multiplicity one in its action on the Lie algebra
$Lie U^+(\C)\simeq \C ^2\otimes \C^{2n-4}\oplus triv$.  By Proposition
\ref{multone}, Theorem 1 holds. \\

{\it  Case   4.   $K=\Q$  and  $D\otimes  \R\neq   {\bf  H}$}.   Then,
$SL_1(D\otimes  \R)$ is  non-compact and  semi-simple. Now,  the group
$SL_1(D)\times    SL_2/K$    is    embedded   in    $SU(J,D)$    where
$J=\begin{pmatrix} 0 & 1\\ 1  & 0\end{pmatrix}$ is the hyperbolic form
in  two variables.   Therefore, the  real rank  of $SU(J,D)$  is $\geq
2$. \\

Write  $h_{n-2}=\lambda  _1\oplus   h_{n-3}$  for  some  $\lambda  \in
D^{\sigma}-\setminus  \{0\}$.  After  a  scaling, we  may assume  that
$\lambda =1$. Consider the group $G_1=SU_3(J\oplus 1, D)$.  Let $P_1$,
$U_1$  be  the intersections  of  $P$ and  $U$  with  $G_1$. They  are
respectively a  parabolic subgroup and  its unipotent radical  in $H$.
By the last paragraph, it follows that $M_0\simeq SL_1(D)$. \\

Now, the Tits  diagram of $G_1$ is that  of $^2A_3\simeq SU(1,3)$ over
$\Q$, where $SU(1,3)$ actually  denotes the $K$-rank one group $SU(B)$
with $B$ a hermitian form in four variables over a quadratic extension
$E$ of  $\Q$, such that the  maximal isotropic subspaces  of $E^4$ for
the form $B$ are one dimensional. Thus, $G_1$ is as in Cases 3 or 4 of
subsection (5.2.1).  In  Case 3 of (5.2.1), it is easy  to see (and is
observed there) that $M_0$ is  not semi-simple. Therefore, only Case 4
of (5.2.1)  applies. In this case  (see (5.2.1), Case 4),  there is an
embedded  $H=SO(1,3)$ in this  $G_1\simeq SU(1,3)$  of real  rank two.
Choose the unipotent element $u_0\in H\cap U_1(\Z)\subset U^+(\Z)$ and
$v_0\in (U_1)_{-2\alpha}(\Z)=U_{2\alpha}(\Z)$ (the last equality holds
since the  space $\fg _{2\alpha}$  is one dimensional) and  an element
$\theta  \in M_0\cap  H$  as  in (5.2.1),  Case  4. Set  $V^+(r)=H\cap
U^+(r\Z)U_{2\alpha}(r\Z)$.  Then, by  the argument of section (5.2.1),
Case  4, $V^+(r)$ is  contained in  the two-generated  group $<(\theta
)^r, (u_0v_0)^r>$.  Let $\gamma \in G(\Z)$ be in general position with
respect to  $u_0v_0$ and  $\theta$.  By Lemma  \ref{zariskidense}, for
each  $r$,  the  group  $\Gamma=<(u_0v_0)^r,\theta ^r,\gamma  ^r>$  is
Zariski  dense.  To  prove Theorem  1, it  is sufficient  (by  the now
familiar  arguments)  to prove  that  $\Gamma  $  is arithmetic.   Let
$V^-(r)$ denote the $w$ conjugate of $V^+(r)$. \\

Pick a  generic element $g=umwv\in  \Gamma $.  $\Gamma $  contains the
group $<^g(V^+), V^+>\supset  ^u<^m(V^-(r')), V^+(r')>$ for some $r'$.
Thus,   $\Gamma  $   contains  the   subgroup  $^u<U_{-2\alpha}(r'\Z),
U_H(r\Z)>$ where $U_H=U^+\cap  H$; it is proved in  Case 4 of (5.2.1),
that  the   group  $<U_{-2\alpha}(r'\Z),U_H(r'\Z)>$  contains  $\theta
^{r''\Z}$  for some  $r''$. Therefore,  $\Gamma $  contains $^u(\theta
^{r''\Z})$. By Proposition \ref{technical}, $\Gamma $ is arithmetic.

{\it Case 5. $K=\Q$ and $D\otimes \R={\bf H}$}.

If, as before,  $J=\begin{pmatrix} 0 & 1\\ 1  & 0\end{pmatrix}$ is the
hyperbolic form in two variables  over the division algebra $D$, then,
$SU(J,D)(\R)=\{g\in  SL_2({\bf  H}):  g^{\sigma}Jg=J\}$ has  $\R$-rank
$1$. Recall  that $h=J\oplus h_{n-2}$. Since  $R$-rank ($SU(h)$) $\geq
2$, we  must have $\R$-rank ($SU(h_{n-2})$) $\geq  1$. Hence $h_{n-2}$
represents a zero over $\R$, and therefore, $n-2\geq 2$.  \\

If  $n-2\geq   3$,  then  write   $h_{n-2}=h_2\oplus  h_{n-4}$.   Now,
$G_0=SU(J\oplus  h_2,D)$  is  an  absolutely simple  $\Q$-subgroup  of
$G$. We  will show that $\Q$-rank  ($G_0$) $\geq 2$,  which will prove
the same for $G$, and contradicts our assumption that $\Q$-rank of $G$
is one.\\

The group  $G_0$ is of type  $D_4$, with $\Q$-rank one.   Thus, In the
diagram of $G_0$, there is one circled root. If the anisotropic kernel
$M'$  is  $\Q$-simple, then,  $M'(\C)\supset  SL_2^3$, and  therefore,
$<U^+_{G_0}(\Z), S_0(\Z)>$  (with $S_0$ a suitable torus  in $M'$), is
two generated: say by $u_+$  and $\theta $.  By considering an element
$\gamma \in G(O_K)$, in general  position it follows that $<\gamma ^r,
\theta  ^r, u_+^r>$  is  Zariski  dense.  Write  $V^+$  for the  group
generated by $\theta  ^r$ and $u_+^r$, and $V^-$  for its conjugate by
$w$. Since  $G_0$ contains the  $2\alpha$ root group  $U_{2\alpha}$ it
follows  that $V^+$ is  normalised by  the unipotent  arithmetic group
$U^+(r\Z)$.   Consequently, given $g=umwv  \in \Gamma  \cap U^+MwU^+$,
$\Gamma  $ contains  the  group $<^g(V^+),V^+>=^u<^m(V^-),V^+>$.   The
latter  contains  $\Delta=^u<U_{-2\alpha}(r\Z),U^+_{G_0}(r\Z)>$. Since
$G_0$ is of higher real rank (and of $\Q$-rank one), any Zariski dense
subgroup  of $G_0(\Z)$ intersecting  $U^+_{G_0}(\Z)$ in  an arithmetic
group is of finite index in $G_0(\Z)$ by \cite{V}. Therefore, $\Delta$
is   of  finite  index   in  $G_0(\Z)$   and  hence   $\Gamma  \supset
^u(S_0(r'\Z))$ for  some integer  $r'$.  Now, non-trivial  elements of
$S_0(r'\Z)$, act  by eigenvalues $\neq  1$ on the $\alpha$  root space
$\fg  _\alpha$.   An argument  similar  to  the  proof of  Proposition
\ref{technical} shows  that the Zariski  closure $\frak v$  of $\Gamma
\cap U^+$  has Lie algebra  which contains $\fg _\alpha$.   The latter
{\it generates}  $\fu$. Therefore,  $\frak v=\fu$ and  $\Gamma \supset
U^+(r''\Z)$  for  some  $r''$.    Thus,  by  \cite{V},  $\Gamma  $  is
arithmetic, and Theorem 1 holds.\\

[2]. If the  anisotropic kernel is not $\Q$-simple,  then, there is at
least  one  simple root  connected  to  the  above circled  root,  and
together,  they  generate a  group  $G_1$  isomorphic  to $SL_3$  over
$\C$. Over $\R$, $G_1$ cannot be  outer type $SL_3$, since one root is
already circled over  $\Q$ (in outer type $A_2$,  {\it two} roots over
$\R$,  are  circled  together).    Therefore,  $G_1$  is  $SL_3$  over
$\R$. Hence, over $\Q$, $G_1$ can  only be $SU(2,1)$ with respect to a
real  quadratic  extension. Then  again,  the group  $<U^+_{G_1}(r\Z),
\theta  ^{r\Z}>$ is  virtually  two  generated (for  any  $r$), and  a
general  position argument  as in  the previous  paragraph  shows that
Theorem 1 holds in this case too. \\

\newpage

\section{Exceptional groups of rank one} 

\subsection{Groups of type $^3D_4$ and $^6D_4$} The only $K$-rank one
groups   (according  to   \cite{T2},  p.58)   are   $^3D_{4,1}^9$  and
$^6D_{4,1}^9$. The simple root that  is connected to all the others is
circled.   The  anisotropic  kernel  $M_1$  is,  over  $\overline  K$,
$SL_2^3$.   Moreover,  the  Galois  group of  $\overline  {K}/K$  acts
transitively on  the roots connected  to this simple root.   Thus, the
anisotropic kernel  is an {\it  inner twist} of the  quasi-split group
$M'=R_{E/K}(SL_2)$ with  $E/K$ either cubic  ($^3D_{4,1}^9$) or sextic
($^6D_{4,1}^9$).   $M'$is   $K$-simple  whence  any   inner  twist  is
$K$-simple (inner twist of a product is a product of inner twists). \\

Now, $G$ being an inner twist of the quasi-split group ${\mathcal G}$,
is given by  an element of the Galois  cohomology set $H^1(K,{\mathcal
G})$.   However,   this  element  is  in  the   image  of  $H^1(K,M')$
(Proposition 4 (ii) of \cite{T2}). Hence $G$ contains the $K$-subgroup
$M_1$ (inner twist of $M'$), whence $M_1$ is $K$-simple. \\

Since $\R$-rank  ($M_1(K\otimes \R)$) $\geq 1$ (it  follows by looking
at the Tits  diagrams, that $G(K_v)$ has $K_v$-rank  $\geq 2$ for each
archimedean place  $v$ of $K$, because  these forms do  not occur over
real  or  complex numbers),  it  follows  that  $M_1(K\otimes \R)$  is
non-compact semi-simple.   Hence the Zariski closure  of $M_1(O_K)$ is
$M_1$.    Now,   by   \cite{L},    \cite{Sh},   as   a   module   over
$M_1(\C)=SL_2(\C)^3$,  $Lie U^+=St\otimes  St\otimes  St$.  Therefore,
the torus of $M_1(\C)$ given by the product of diagonal tori in $SL_2$
acts by multiplicity one  on $Lie U^+$.  By Proposition \ref{multone},
Theorem 1 follows.

\subsection{Groups of type $E_6$} 

{\it Case 1. There are no inner type groups of rank one}.

{\it Case  2.  $G=^2E_{6,1}^{35}$}.   The anisotropic kernel  $M_1$ is
$K$-simple  (since  its  Tits  diagram  is  connected).   It  is  also
non-compact  at infinity, since  $\R$rank of  any non-compact  form of
$E_6$  over  $K_v$  has   $K_v$-rank  $\geq  2$  for  any  archimedean
completion  of  $K$.   Hence  $M_0\supset  M_1$.   As  a  module  over
$M_1(\C)=SL_6(\C)$,  $Lie  U^+$ is  $\wedge  ^3 (\C^{10})$  (\cite{L},
\cite{Sh}),  and  is  multiplicity  free  for the  diagonal  torus  in
$SL_6$. This completes the proof.\\

{\it  Case   3.   $G=^2E_{6,1}^{29}$}.   The   anisotropic  kernel  is
non-compact  at  infinity  for   the  same  reason  as  above.   Hence
$M_1\subset M_0$,  with $M_1=SO(8)$. As  an $M_1(\C)=SO(8,\C)$ module,
the space $Lie U^+$ is  (by p. 568, $^2E_6$-3 of \cite{Sh}), $St\oplus
\delta _3\oplus \delta _5$ where $St$, $\delta _3$ and $\delta _5$ are
respectively, the  standard, and the  two distinct spin  modules. With
respect  to   the  maximal  torus  of  $SO(8,\C)$,   the  weights  are
$x_1,\cdots x_4, \frac{\epsilon _1x_1+\cdots+ \epsilon _4x_4}{2}$ with
$\epsilon _i=\pm 1$, each  occurring with multiplicity one. Therefore,
Theorem 1 follows from Proposition \ref{multone}.\\

\subsection{Groups of type $E_7$ or $E_8$ or $G_2$}
 
There are no $K$-rank one forms over number fields.

\subsection{Groups of type $F_4$}

The $K$-rank one form is  $F_{4,1}^{21}$. This is the only exceptional
group  which can  have rank  one over  some archimedean  completion of
$K$. \\

{\it Case  1. $K$  is not  totally real or  $K=\Q$ or  the anisotropic
kernel  is non-compact  at  infinity}.  Then,  the anisotropic  kernel
$M_1$ is a  form of $SO(7)$.  Over $\C$ this  is non-compact.  In case
$K=\Q$ again, this is non-compact over  $\R$ since $G$ is of real rank
$\geq 2$.   If $K\neq \Q$ is  totally real, then  by {\it assumption},
$M_1$  is non-compact  at  infinity.  Thus,  $Lie U^+=St\oplus  \wedge
^3(\C^7)$ (\cite{L}, (xxii), p.52)  is multiplicity free for the torus
of $SO(7)$. \\

{\it Case 2.  $K$ totally  real, and the anisotropic kernel is compact
at infinity}. let $\fg _{2\alpha}$  be the $2\alpha$ root space. Then,
the   subgroup  $G_1$   of   $G$  with   Lie   algebra  $\fg   _1=<\fg
_{-2\alpha},\fg_{2\alpha}>$    must   be    locally    isomorphic   to
$SO(1,8)$. For,  $\fg _1$  has real  rank one (since  $\fg $  has), is
semi-simple, and its obvious  parabolic subgroup has abelian unipotent
radical.   Therefore,  it  can  only  be  $SO(1,k)$.   Since  $dim(\fg
_{2\alpha})=7$,  it  follows  that  $k-1=7$  i.e.   $k=8$.   Now,  the
anisotropic factor  $SO(7)$ of $G_1=SO(1,8)$ is  an anisotropic factor
of $F_{4,1}^{21}$ as well. \\

Fix  $u_+\in  U_{2\alpha}(O_K)$,   and  $\theta  \in  {\bf  G}_m(O_K)$
suitably chosen (as in Lemma  \ref{exist}). Fix $\gamma \in G(O_K)$ in
general  position   with  respect  to  $u_+,\theta   $  (  Proposition
\ref{zariskidense}).   For  each  $r$,  write  $\Gamma  =<u_+^r,\theta
^r,\gamma ^r>$.   Then, 1) $\Gamma  $ is Zariski dense  in $G(K\otimes
\C)$        (        Proposition       \ref{zariskidense}).         2)
$V^+=V^+(r')==u_+^{r'O_K}\subset \Gamma  $ for some  integer $r'$. Put
$V^-=V^-(r')$ for  the $w$ conjugate  of $V^+(r')$. 3) If  $g=umwv \in
\Gamma  $  is   generic,  then  $\Gamma  \supset  <^g(V^+),V^+>\supset
^u<^m(V^-),V^+>$.   By using the  result proved  for $SO(1,8)$  (it is
important  to note that  $K\neq \Q$  is totally  real, and  that $m\in
SO(7)\subset SO(1,8)$ to  apply this result), we see  that $^u (\theta
^{r''\Z})\subset \Gamma $ for  some $r''\neq 0$.  Then, by Proposition
\ref{technical}, $\Gamma $ is arithmetic.\\

\end{document}